\documentclass[12pt,leqno]{amsart}

\usepackage{amssymb}
\usepackage{amsthm}
\usepackage{amsmath}
\usepackage{amscd}
\usepackage[mathscr]{eucal}
\usepackage{verbatim}
\usepackage[all]{xy}

\numberwithin{equation}{section}

\theoremstyle{plain}
\newtheorem{theorem}{Theorem}[section]
\newtheorem*{theorem*}{Theorem}
\newtheorem{corollary}[theorem]{Corollary}
\newtheorem*{corollary*}{Corollary}
\newtheorem{lemma}[theorem]{Lemma}
\newtheorem{proposition}[theorem]{Proposition}

\theoremstyle{definition}

\newtheorem{remark}[theorem]{Remark}
\newtheorem{example}[theorem]{Example}

\theoremstyle{remark}

\newcommand{\A}{\mathbb{A}}
\newcommand{\R}{\mathbb{R}}
\newcommand{\Q}{\mathbb{Q}}
\newcommand{\Z}{\mathbb{Z}}

\newcommand{\C}{\mathbb{C}}
\newcommand{\h}{\mathbb{H}}
\renewcommand{\H}{\mathbb{H}}

\newcommand{\G}{\Gamma}
\newcommand{\g}{\gamma}

\newcommand{\La}{\Lambda}
\newcommand{\x}{\mathbf{x}}

\newcommand{\back}{\backslash}

\newcommand{\zxz}[4]{\begin{pmatrix} #1 & #2 \\ #3 & #4 \end{pmatrix}}

\newcommand{\kzxz}[4]{\left(\begin{smallmatrix} #1 & #2 \\ #3 & #4\end{smallmatrix}\right) }

\newcommand{\re}{\operatorname{Re}}

\newcommand{\calO}{\mathcal{O}}

\newcommand{\calS}{\mathcal{S}}

\newcommand{\calZ}{\mathcal{Z}}

\newcommand{\frake}{\mathfrak e}

\newcommand{\bs}{\backslash}

\newcommand{\Ad}{\operatorname{Ad}}
\newcommand{\Lie}{\operatorname{Lie}}

\newcommand{\vol}{\operatorname{vol}}

\newcommand{\Span}{\operatorname{span}}

\newcommand{\Sl}{\operatorname{SL}}

\newcommand{\SL}{\operatorname{SL}}
\newcommand{\GL}{\operatorname{GL}}

\newcommand{\Sp}{\operatorname{Sp}}

\newcommand{\Mat}{\operatorname{Mat}}
\newcommand{\Mp}{\operatorname{Mp}}
\newcommand{\Orth}{\operatorname{O}}
\newcommand{\Uni}{\operatorname{U}}

\newcommand{\Sym}{\operatorname{Sym}}

\newcommand{\SO}{\operatorname{SO}}

\newcommand{\Pic}{\operatorname{Pic}}

\newcommand{\dv}{\operatorname{div}}

\begin{document}
\title[On the Kudla-Millson lift]{On the injectivity of the Kudla-Millson lift and surjectivity of the Borcherds lift}
\author[Jan H.~Bruinier and Jens Funke]{Jan Hendrik Bruinier and Jens Funke*}
\thanks{* Partially supported by NSF grant  DMS-0305448}
\address{Mathematisches Institut, Universit\"at zu K\"oln, Weyertal 86--90, D-50931 K\"oln, Germany}
\email{bruinier@math.uni-koeln.de }
\address{Department of Mathematical Sciences, New Mexico State University, P.O.~Box 30001, 3MB, Las Cruces, NM 88003, USA}
\email{jfunke@nmsu.edu}
\date{June 5, 2006}

\begin{abstract}
We consider the Kudla-Millson lift from elliptic modular forms of
weight $(p+q)/2$ to closed $q$-forms on locally symmetric
spaces corresponding to the orthogonal group $\Orth(p,q)$.  We study
the $L^2$-norm of the lift following the Rallis inner product
formula. We compute the contribution at the Archimedian place. For
locally symmetric spaces associated to even unimodular lattices, we
obtain an explicit formula for the $L^2$-norm of the lift, which
often implies that the lift is injective.  For $\Orth(p,2)$ we
discuss how such injectivity results imply the surjectivity of the
Borcherds lift.
%In particular, we revisit part of the material of
%\cite{Br2} from the viewpoint of \cite{BF}.
\end{abstract}

\maketitle

\section{Introduction}

In previous work \cite{BF}, we studied the Kudla-Millson theta lift
(see e.g.~\cite{KM90}) and Borcherds' singular theta lift
(e.g.~\cite{Bo2,Br2}) and established a duality statement between
these two lifts. Both of these lifts have played
a significant role in the study of certain cycles in locally
symmetric spaces and Shimura varieties of
orthogonal type. In this paper, we study the injectivity of the
Kudla-Millson theta lift, and revisit part of the material of
\cite{Br2} from the viewpoint of \cite{BF},
% use the results of \cite{Br2} and \cite{BF}
to obtain surjectivity results for the Borcherds  lift. Moreover, we
provide evidence for the following principle: The vanishing of the
standard $L$-function of a cusp form of weight $1+p/2$ at $s_0=p/2$
corresponds to the existence of a certain "exceptional automorphic
product" on $\Orth(p,2)$ (see Theorem~\ref{borcherdsexceptional}).

We now describe the content of this paper in more detail. We begin by
recalling the Kudla-Millson lift in a setting which is convenient
for the application to the Borcherds lift.
Let $(V,Q)$ be a non-degenerate rational quadratic space of
signature $(p,q)$. We write $(\cdot,\cdot)$ for the bilinear form
corresponding to the quadratic form $Q$. We write $r$ for the Witt
index of $V$, i.e., the dimension of a rational
maximal isotropic subspace.
%over $\Q$.
Throughout we assume that the dimension $m=p+q$ of $V$ is
even.  Let $H=\Orth(V)$ be the orthogonal group of $V$.
%and write $H(\R)_0$ for the connected component of the identity of
%$H(\R)$.
We let $D$ be the bounded domain associated to $H(\R)$, which we realize as
the Grassmannian of oriented negative $q$-planes in $V(\R)$.
%\[
%D=\{z\subset V(\R);\quad\text{$\dim z=q$ and  $Q|_z<0$}\}.
%\]

Let $L\subset V$ be an even lattice of level $N$, and write $L^{\#}$ for
the dual lattice. The quadratic form on $L$ induces a non-degenerate
$\Q/\Z$-valued quadratic form on the discriminant group $L^{\#}/L$.
Recall that the Weil representation $\rho_L$ of the quadratic module
$(L^{\#}/L,Q)$ is a unitary representation of $\Sl_2(\Z)$ on the group
ring $\C[L^{\#}/L]$, which can be defined as follows \cite{Bo2},
\cite{Br2}.  If $(\frake_\gamma)_{\gamma\in L^{\#}/L}$ denotes the
standard basis of $\C[L^{\#}/L]$, then $\rho_L$ is given by the action
of the generators $T=\kzxz{1}{1}{0}{1}$ and $S=\kzxz{0}{-1}{1}{0}$
of $\Sl_2(\Z)$ by
\begin{align*}
%\label{eq:weilt}
\rho_L(T)(\frake_\gamma)&=e(\gamma^2/2)\frake_\gamma,\\
%\label{eq:weils}
\rho_L(S)(\frake_\gamma)&=\frac{e(-(p-q)/8)}{\sqrt{|L^{\#}/L|}}
\sum_{\delta\in L^{\#}/L} e(-(\gamma,\delta)) \frake_\delta,
\end{align*}
where $e(w):=e^{2\pi i w}$. This representation factors through the
group $\Sl_2(\Z/N\Z)$.

Let $\Gamma\subset \Orth(L)$ be a torsion-free subgroup of finite
index which acts trivially on $L^{\#}/L$. Then
\[
X=\Gamma\bs D
\]
is a real analytic manifold. For $x \in L^{\#}$ with $Q(x) >0$, we let
\begin{equation*}
D_x = \{ z \in D;\; z \perp x\}.
\end{equation*}
Note that $D_x$ is a subsymmetric space attached to the orthogonal
group $H_x$, the stabilizer of $x$ in $H$. Put $\Gamma_x = \Gamma \cap
H_x$. The quotient
\begin{equation*}
Z(x) = \G_x \back D_x \longrightarrow X
\end{equation*}
defines a (in general relative) cycle in $X$. For $h \in
L^{\#}/L$ and $n\in \Q$, the group $\G$ acts on $L_{h,n} = \{x \in L +h;\; Q(x)
= n \}$ with finitely many orbits, and we define the
composite cycle
\begin{equation*}
Z(h,n) = \sum_{ x \in \G \back L_{h,n}} Z(x).
\end{equation*}
%Occasionally, we will identify  $Z(h,n)$ with its preimage in $D$.

Kudla and Millson constructed Poincar\'e dual forms for such cycles
by means of the Weil representation, see e.g.~\cite{KM90}.
They constructed a Schwartz form
$\varphi_{KM}\in[\calS(V(\R))\otimes \mathcal{Z}^q(D)]^{H(\R)}$ on
$V(\R)$ taking values in $\mathcal{Z}^q(D)$, the closed differential
$q$-forms on $D$.
%with the property that $\varphi_{KM}(x,z)e^{\pi
%(x,x)}$ is a Poincar\'e dual form for the cycle $D_x$.
Let $\omega_\infty$ be the Schr\"odinger model of the Weil
representation of $\Sl_2(\R)$ acting on the space of Schwartz
functions $\calS(V(\R))$, associated to the standard additive
character. We obtain a $\C[L^{\#}/L]$-valued theta function on the
upper half plane $\H$ by putting
\begin{align*}
\Theta(\tau,z,\varphi_{KM}) = v^{-m/4}\sum_{h\in L^{\#}/L}\sum_{x \in
L+h} (\omega_\infty(g_\tau)\varphi_{KM})(x,z)\frake_h.
\end{align*}
Here $\tau=u+iv\in \H$ and $g_\tau= \kzxz{1}{u}{0}{1}
\kzxz{\sqrt{v}}{0}{0}{\sqrt{v}^{-1}} \in \Sl_2(\R)$ is the standard
element moving the base point $i \in \h$ to $\tau$.
In the variable $\tau$, this theta function transforms
as a (non-holomorphic) modular form of weight $\kappa=m/2$ for
$\Sl_2(\Z)$ of type $\rho_L$.
In the variable $z$, it defines a closed $q$-form on $X$.
Kudla and Millson showed that the Fourier coefficient
at $e^{2\pi i n \tau}\frake_h$
is a Poincar\'e dual form for the cycle $Z(h,n)$.

Let $S_{\kappa,L}$ denote the space of $\C[L^{\#}/L]$-valued cusp forms
of weight $\kappa$ and type $\rho_L$ for the group  $\Sl_2(\Z)$. We
define a lifting $\Lambda:S_{\kappa,L}\to \calZ^q(X)$ by the theta
integral
\begin{align}
\label{introkm}
f\mapsto \Lambda(f)=\int_{\Sl_2(\Z)\bs \H} \langle
f(\tau), \Theta(\tau,z,\varphi_{KM})\rangle \,\frac{du\,dv}{v^2},
\end{align}
where $\langle \cdot, \cdot \rangle$ denotes the standard scalar
product on $\C[L^{\#}/L]$.
%One can in fact show that $\Lambda(f)$ is a harmonic form.
%We have chosen the present formulation of the lift
%since it fits well the application to the Borcherds lift below.

In the present paper, we consider the question whether $\Lambda$ is
injective.  We compute the $L^2$-norm of the differential form
$\Lambda(f)$ in the sense of Riemann geometry by means of the Rallis
inner product formula \cite{Ra}.
%We let $\Sp(2)$ be the symplectic
%group of rank $2$. Then,
First, using the see-saw
\[
\xymatrix{
\Sp(2) \ar@{-}[d] \ar@{-}[dr]& \Orth(V)\times \Orth(V)\ar@{-}[d]\\
\SL_2\times \SL_2\ar@{-}[ur]& \Orth(V)
}
\]
and the Siegel-Weil formula (see e.g.~\cite{KR1}, \cite{KR2},
\cite{KR3}, \cite{Weil}), the inner product can be expressed as a
convolution integral of $f$ against the restriction of a genus $2$
Eisenstein series to the diagonal (see Proposition \ref{prop:l21}).

Such convolution integrals can be evaluated by means of the doubling
method, see e.g.~\cite{Boe}, \cite{Ga}, \cite{PS-R}, \cite{Ra}. If $f$
is a Hecke eigenform of level $N$, one obtains a special value of the
partial standard $L$-function of $f$ (where the Euler
factors corresponding to the primes dividing the level $N$ and $\infty$
are omitted)
times a product of ``bad''
local factors corresponding to the primes dividing $N$ and
$\infty$.  If $m>4$, then, by the Euler product expansion, the special
value of the partial standard $L$-function is non-zero.  Therefore the lift
$\Lambda(f)$ vanishes precisely if at least one of the ``bad'' local
factors vanishes.  By the analysis of the present paper we determine
the local factor at infinity.
%This is done in Section
%\ref{Schwartz} by a computation in the Fock model of the Weil
%representation (see in particular Proposition \ref{decomp} and
%Proposition \ref{keyprop}).

In the special case where $L$ is even and unimodular, the level of $L$
is $N=1$, so that $\infty$ is the only ``bad'' place. The space
$S_{\kappa,L}$ is equal to the space $S_\kappa(\Gamma(1))$ of scalar
valued cusp forms of weight $\kappa$ for $\Gamma(1)=\Sl_2(\Z)$.  We
obtain the following explicit formula for the $L^2$-norm of the lift
(see Theorem \ref{thm:main}):

\begin{theorem}
\label{intromain}
  Assume that $m>3+r$, where $r$ is the Witt index of $V$.  Let $f\in
  S_\kappa(\Gamma(1))$ be a Hecke eigenform, and write
  $\|f\|^2_2$ for its Petersson norm, and $D_f(s)$ for its
  standard $L$-function.  Then $\Lambda(f)$ is square integrable and
\begin{align*}
\frac{\|\Lambda(f)\|_2^2}{\|f\|^2_2} = C\cdot
\frac{D_f(m/2-1)}{\zeta(m/2)\zeta(m-2)},
\end{align*}
where $C=C(p,q)$ is an explicit real constant, which does not depend on $f$.%(which is computed in the body of the paper).
The constant $C$ vanishes if and only if $p=1$.
\end{theorem}

%  By the same argument it is easily seen that
%  $(\Lambda(f),\Lambda(g))=0$ for two different normalized Hecke
%  eigenforms $f$ and $g$.

\begin{corollary}
Assume that $m>\max(4,3+r)$ and that $L$ is even unimodular. When
$p\neq 1$, the theta lift $\Lambda$ is injective. When $p=1$, the
lift vanishes identically.
\end{corollary}

It would be interesting to compute the bad local factors at finite
primes (or at least to show their non-vanishing) as well. However,
in our setting, this requires first a suitable Hecke theory for
vector valued modular forms in $S_{\kappa,L}$, which is not yet
available. It seems conceivable that one could prove more general
injectivity results along these lines. For the relationship between
the vector-valued modular forms in $S_{\kappa,L}$ and the adelic
language, see \cite{KBforms}.
%The present paper can be viewed as a first step in that direction,
%carrying out the computation at $\infty$.

Note that in this context, J.-S. Li \cite{Li} has used the theta
correspondence and the doubling method for automorphic
representations in great generality to obtain non-vanishing results
for cohomology when passing to a sufficiently large level.

In the body of the paper, we actually consider the generalization of
the Kudla-Millson lift due to Funke and Millson \cite{FM}.
It maps cusp forms in
$S_{\kappa,L}$ to closed differential $q$-forms with values in certain local
coefficient systems.
Moreover, we use an adelic set-up for the theta and Eisenstein
series in question.

\subsection{Surjectivity of the Borcherds lift}

We briefly discuss how the injectivity results on the Kudla-Millson
lift imply surjectivity results for the Borcherds lift. We revisit
part of the material of \cite{Br2} in the light of the adjointness
result of \cite{BF} between the regularized theta lift and the Kudla
Millson lift.
%As a Corollary to the above Theorem we obtain yet
%another proof of the surjectivity of the Borcherds lift for
%unimodular lattices.
We restrict ourselves to the Hermitean case of
signature $(p,2)$ where $X$ is a $p$-dimensional complex algebraic
manifold. The special cycles $Z(h,n)$ are algebraic divisors on $X$,
also called Heegner divisors or rational quadratic divisors.

We say that a meromorphic modular form for $\Gamma$ has a Heegner
divisor, if its divisor on $X$ is a linear combination of the
$Z(h,n)$. A large supply of modular forms with Heegner divisor is
provided by the Borcherds lift, see \cite{Bo1}, \cite{Bo2}. We
briefly recall its construction.

A meromorphic modular form for a congruence subgroup of
$\Sl_2(\Z)$ is called {\em weakly holomorphic}, if its poles are
supported on the cusps. If $k\in \Z$, we write $W_{k,L}$ for the
space of weakly holomorphic modular forms of weight $k$ for
$\SL_2(\Z)$ of type $\rho_L$. Any $f\in W_{k,L}$ has a Fourier
expansion of the form
\[
f(\tau)= \sum_{h\in L^{\#}/L} \sum_{n\in \Z+Q(h)} c(h,n)
e(n\tau)\frake_h,
\]
where only finitely many coefficients $c(h,n)$ with $n<0$ are
non-zero. We write $V^-$ for the quadratic space $(V,-Q)$ of
signature $(2,p)$ and $L^-$ for the lattice $(L,-Q)$ in $V^-$.

\begin{theorem}[Borcherds \cite{Bo2}, Theorem 13.3]
\label{introborcherds}
Let $f\in W_{1-p/2,L^-}$ be a weakly
holomorphic modular form with Fourier coefficients $c(h,n)$. Assume
that $c(h,n)\in \Z$ for $n<0$. Then there exists a meromorphic
modular form $\Psi(z,f)$ for $\Gamma$ (with some multiplier system
of finite order) such that:
\begin{enumerate}
\item[(i)]
The weight of $\Psi$ is equal to $c(0,0)/2$.
\item[(ii)]
The divisor $Z(f)$ of $\Psi$ is determined by the principal part of
$f$ at the cusp $\infty$. It equals
\[Z(f)=\sum_{h\in L^{\#}/L}\sum_{n<0} c(h,n)Z(h,n).\]
\item[(iii)]
In a neighborhood of a cusp of $\Gamma$
%(given by an isotropic line in $V^-$)
the function $\Psi$ has an infinite product expansion
analogous to the Dedekind eta function, see \cite{Bo2}, Theorem 13.3 (5.).
\end{enumerate}
\end{theorem}

The proof of this result uses a regularized theta lift. Let
$\varphi_0^{p,2}\in \calS(V(\R))$ be the Gaussian for signature $(p,2)$. The corresponding
Siegel theta function
\begin{align*}
\Theta(\tau,z,\varphi_{0}) = v \sum_{h\in L^{\#}/L}\sum_{x \in L+h}
(\omega_\infty(g_\tau)\varphi_{0})(x,z)\frake_h
\end{align*}
transforms like a non-holomorphic modular form of weight $p/2-1$ of
type $\rho_L$ in the variable $\tau$. Hence the theta integral
\begin{equation}\label{defbor}
\Phi(z,f) = \int_{\Gamma(1) \back \H} \langle f(\tau),
\overline{\Theta(\tau,z,\varphi_0^{p,2})} \rangle \,d \mu
\end{equation}
formally defines a $\Gamma$-invariant function on $D$. Because of
the singularities of $f$ at the cusps, the integral diverges.
However, Harvey and Moore discovered that it can be regularized
essentially by viewing it as the limit $T\to \infty$ of the integral
over the standard fundamental domain truncated at $\Im(\tau)=T$, see
\cite{Bo2}, \cite{HM}. It turns out that $\Phi(z,f)$ defines a
smooth function on $X\setminus Z(f)$ which  has a logarithmic
singularity along $Z(f)$. Moreover,
\begin{align*}
%\label{reg-theta}
\Phi(z,f)=-2\log\| \Psi(z,f)\|_{\text{Pet}} +
\text{constant},
\end{align*}
where $\|\cdot\|_{\text{Pet}}$ denotes the Petersson metric on the
line bundle of modular forms of weight $c(0,0)/2$ over $X$. From
this identity, the claimed properties of $\Psi(z,f)$ can be derived.

Modular forms for the group $\Gamma\subset \Orth(L)$ arising via
this lift are called automorphic products or Borcherds products. By
(ii) they have a Heegner divisor.

Here we consider the question whether the Borcherds lift is
surjective. More precisely we ask whether every meromorphic modular
form for $\Gamma$ with Heegner divisor is the lift $\Psi(z,f)$ of a
weakly holomorphic form  $f\in W_{1-p/2,L^-}$?

An affirmative answer to this question was given in \cite{Br2} in
the special case that the lattice $L$ splits two hyperbolic planes over $\Z$. In the
(more restrictive) case that $L$ is unimodular, a different proof was
given in \cite{BrFr} using local Borcherds products and a theorem of
Waldspurger on theta series with harmonic polynomials \cite{Wald}.

The approach of \cite{Br2} was to first simplify the problem and to
consider the regularized theta lift for a larger space of ``input''
modular forms. Namely, we let $H_{k,L}$ be the space of \emph{weak
Maass forms} of weight $k$ and type $\rho_L$.  This space consists
of the smooth  functions $f:\h\to \C[L^{\#}/L]$ that transform with
$\rho_L$ in weight $k$ under $\Sl_2(\Z)$, are annihilated by the
weight $k$ Laplacian, and satisfy $f(\tau) = O(e^{Cv})$ as $\tau= u+
iv \to i\infty$ for some constant $C>0$ (see \cite{BF} Section 3).

Any $f\in H_{k,L}$ has a Fourier expansion of the form
\begin{align}
\label{weakmaass}
f(\tau)&= \sum_{h\in L^{\#}/L}\sum_{n\in \Q} c^+(h,n) e(n\tau)\frake_h\\
\nonumber &\phantom{=}{}+ \sum_{h\in L^{\#}/L} c^-(h,0) v^{1-k}\frake_h
+\sum_{\substack{n\in \Q\\ n\neq 0}} c^-(h,n) H(2\pi nv)
e(nu)\frake_h,
\end{align}
where only finitely many of the coefficients $c^+(h,n)$
(respectively $c^-(h,n)$) with negative (respectively positive)
index $n$ are non-zero. The function $H(w)$ is a Whittaker type
function.

For $f \in H_{k,L}$, put $\xi_k(f) = R_{-k}(v^k\bar{f})$, where
$R_{-k}$ is the standard raising operator for modular forms of
weight $-k$. This defines  an antilinear map $\xi_k: H_{k,L} \to
W_{2-k,L^-}$ to the space of weakly holomorphic modular forms in
weight $2-k$. It is easily checked that $W_{k,L}$ is the kernel of
$\xi_k$. According to \cite{BF}, Corollary 3.8, the sequence
\[
\xymatrix{ 0\ar[r]& W_{k,L} \ar[r]& H_{k,L} \ar[r]^{\xi_k}&
W_{2-k,L^-} \ar[r] & 0 }
\]
is exact. We let $H_{k,L}^+$ be the preimage under $\xi_k$ of the space of
cusp forms $S_{2-k,L^-}$ of weight $2-k$ with type $\rho_{L^-}$.
Hence we have the exact sequence
\[
\xymatrix{ 0\ar[r]& W_{k,L} \ar[r]& H_{k,L}^+ \ar[r]^{\xi_k}&
S_{2-k,L^-} \ar[r] & 0 }.
\]
The space $H_{k,L}^+$ can also be characterized as the subspace of
those $f\in H_{k,L}$ whose Fourier coefficients $c^-(h,n)$ with
non-negative index $n$ vanish. This implies that
\[
f(\tau)= \sum_{h\in L^{\#}/L} \sum_{n<0} c^+(h,n) e(n\tau)\frake_h +
O(1), \qquad \Im(\tau)\to \infty,
\]
i.e., the singularity at $\infty$, called the {\em principal part} of
$f$, looks like the singularity of a weakly holomorphic form.

For $f\in H_{1-p/2,L^-}$, we can define the regularized theta lift
$\Phi(z,f)$ as in \eqref{defbor}, see \cite{Br2}, \cite{BF}.  This
generalized lift is related to the Kudla-Millson lift $\Lambda$
defined in \eqref{introkm} in the following way (see
\cite{BF}, Theorem 6.1).

\begin{theorem}
\label{introbf}
Let $f\in H_{1-p/2,L^-}^+$ and denote its Fourier expansion as in \eqref{weakmaass}.
The $(1,1)$-form $dd^c \Phi(z,f)$ can be continued to a smooth form on
$X$. It satisfies
\[
dd^c \Phi(z,f)= \Lambda(\xi_{1-p/2}(f))(z) + c^+(0,0)\Omega.
\]
Here $\Omega$ denotes the invariant
K\"ahler form on $D$ normalized as in \cite{BF}.
\end{theorem}

On the other hand, the
following ``weak converse theorem'' is proved in \cite{Br2},
Theorem~4.23.

\begin{theorem}
\label{introbr} Assume that $p>r$. Let $F$ be a meromorphic modular
form for the group $\Gamma$ with Heegner divisor
\[
\dv(F)=\sum_h\sum_{n<0} c^+(h,n) Z(h,n)
\]
(where $c^+(h,n)=c^+(-h,n)$ without loss of generality). Then there
is a weak Maass form $f\in H_{1-p/2,L^-}^+$ with principal part $
\sum_h\sum_{n<0} c^+(h,n) e(n\tau)\frake_h$ whose regularized theta
lift satisfies
\begin{align}
\label{wm}
\Phi(z,f)=-2\log\| F\|_{\text{\rm Pet}} +
\text{\rm constant}.
\end{align}
\end{theorem}

Note that the proof in \cite{Br2} is only given in the case that
$p\geq 3$ (where the assumption on the Witt index is automatically
fulfilled). However, the argument extends to the low dimensional
cases. It is  likely that the hypothesis on the Witt index can be
dropped as well, but we have not checked this.

\begin{corollary}
\label{cor1} Assume that $p>r$. Let $F$ be a meromorphic modular
form for the group $\Gamma$ with Heegner divisor as in Theorem
\ref{introbr}. Let $f \in H_{1-p/2,L^-}^+$ be a weak Maass form
whose regularized theta lift satisfies \eqref{wm}. Then
\[
\Lambda(\xi_{1-p/2}(f))=0.
\]
\end{corollary}

\begin{proof}
The assumption on $f$ implies that
\begin{align*}
%\label{weknow}
dd^c \Phi(z,f) = -2dd^c \log\| F\|_{\text{\rm Pet}} =
c^+(0,0)\Omega.
\end{align*}
On the other hand, according to Theorem \ref{introbf}, we have
\[
dd^c \Phi(z,f)= \Lambda(\xi_{1-p/2}(f))(z) + c^+(0,0)\Omega.
\]
If we combine these identities, we obtain the claim.
\end{proof}

\begin{corollary}
\label{cor2}
Assume the hypotheses of Corollary \ref{cor1}.
If $\Lambda$ is injective, then $f$ is weakly holomorphic,
and $F$ is a constant
multiple of the Borcherds lift
$\Psi(z,f)$ of $f$ in the sense of
Theorem \ref{introborcherds}.
\end{corollary}

\begin{proof}
  By Corollary \ref{cor1} we have $\Lambda(\xi_{1-p/2}(f)) =0$. Since
  $\Lambda$ is injective, we find that $\xi_{1-p/2}(f)=0$. But this
  means that $f$ is weakly holomorphic.
\end{proof}

When the lattice $L$ splits two hyperbolic planes over $\Z$,
it was
proved in \cite{Br2} that
$\Lambda$ is injective
by considering the Fourier expansion of the lift.
In Section \ref{sect:5} of
the present paper we show (for even
unimodular lattices) how such injectivity results can be obtained by
the Rallis inner products formula.

We end this section by stating a converse of Corollary \ref{cor1}.
If $r>0$, we let $\ell\in L$ be a primitive isotropic vector, and
let $\ell'\in L^{\#}$ be a vector with $(\ell,\ell')=1$. We let $L_0$ be
the singular lattice $L\cap \ell^\perp$ and let $K$ be the
Lorentzian lattice $L_0/\Z\ell$.

\begin{theorem}
\label{borcherdsexceptional}
  Assume that $p\geq 2$ and $p>r$. Let $f\in H_{k,L^-}^+$ and assume
  that the Fourier coefficients $c^+(h,n)$ ($n<0$) of the principal part
  of $f$ are integral.  If $\xi_{1-p/2}(f)\in \ker(\Lambda)$, then
  there exists a meromorphic modular form $F$ for $\Gamma$ (with some
  multiplier system of finite order) such that:
\begin{enumerate}
\item[(i)]
The weight of $F$ is equal to $c^+(0,0)/2$.
\item[(ii)]
The divisor of $F$ is equal to
\[Z(f)=\sum_{h\in L^{\#}/L}\sum_{n<0} c^+(h,n)Z(h,n).\]
\item[(iii)]
In a neighborhood of a cusp of $\Gamma$,
given by a primitive isotropic vector $\ell\in L$,
the function $F$ has an automorphic product expansion
\[
F(z)=C e((\rho,z)) \prod_{\substack{\lambda\in K'\\
(\lambda,W)>0}} \prod_{\substack{\delta\in L^{\#}/L\\ \delta |
L_0=\lambda}} \big(1-e((\lambda,z))+(\delta,
\ell'))\big)^{c^+(\lambda,Q(\lambda))}.
\]
Here $C$ is a non-zero constant, and we have used the notation of \cite{Bo2}.
\end{enumerate}
\end{theorem}

\begin{proof}
Theorem \ref{introbf} and the fact that
$\Lambda\left(\xi_{1-p/2}(f)\right) =0$ imply that
\[
dd^c\Phi(z,f) = c^+(0,0) \Omega.
\]
(In particular, if $c^+(0,0)=0$, then $f$ is pluriharmonic.) Now we
can argue as in \cite{Br2}, Lemma~3.13 and Theorem~3.16 to prove the
claim.
\end{proof}

We note that the assumption on $r$ and $p$ is needed to guarantee that
the multiplier system of $F$ has finite order. (When $f$
is not weakly holomorphic, we cannot argue with the embedding trick as in \cite{Bo3}, Correction).

If $f$ is weakly holomorphic, then $\xi_{1-p/2}(f)=0$ and the Theorem
reduces to Theorem \ref{introborcherds}. However, if $\Lambda$ is not
injective, and $f$ is a weak Maass form such that $\xi_{1-p/2}(f)$ is a
non-trivial element of the kernel, then Theorem
\ref{borcherdsexceptional} leads to {\em exceptional automorphic
products}. If there are any cases where $\Lambda$ is not injective, it
would be very interesting to construct examples of such exceptional automorphic
products.

\begin{comment}
In view of Theorem \ref{intromain} the results of the present paper
suggest that the existence (respectively non-existence) of exceptional
automorphic products is dictated by the vanishing (respectively
non-vanishing) at $s=m/2-1$ of the standard $L$-function of cusp forms
in $S_{1+p/2,L}$.
\end{comment}

\begin{remark}
If $p\geq 4$, the existence of the meromorphic modular form $F$ with
divisor (ii)
%in Theorem \ref{borcherdsexceptional},
is related to the fact that
$H^1(X,\calO_X)=0$ in this case, which can be proved following the argument of
\cite{Fr} \S3.1. Therefore the Chern class map
$\Pic(X)\to H^2(X,\Z)$ is injective.
\end{remark}

\begin{comment}
\begin{remark}
Note that in view of Theorem \ref{intromain} and the algebraicity of the special values of the standard $L$-function (modulo a period), see e.g.~\cite{Za} Theorem 2, the kernel of $\Lambda$ should have a basis of modular forms with rational Fourier coefficients.
\end{remark}
\end{comment}

\bigskip

We thank S.~B\"ocherer, E.~Freitag, W.~T.~Gan, S.~Kudla, and
J.~Millson for very helpful conversations on the content of this
paper. The second named author also thanks the Max Planck Institut
f\"ur Mathematik in Bonn/Germany for its hospitality during the
summer 2005 where substantial work on this paper was done.

\section{Theta functions and the Siegel-Weil formula}\label{thetaSW}

Let $V$ be a vector space over $\Q$ of dimension $m$ with a
non-degenerate bilinear form $(\,,\,)$. For simplicity we assume
that $m$ is even. We let $H=\Orth(V)$ be the orthogonal group of $V$,
and we let $G=\Sp(n)$ be the symplectic group acting on a symplectic
space of dimension $2n$ over $\Q$. The embedding of $\Uni(n)$ into
$G(\R)$ given by $\mathbf{k}=A + iB \mapsto k=\left(
\begin{smallmatrix}
A & B \\
-B &A
\end{smallmatrix}\right)$ gives rise to a maximal compact subgroup
$K_{\infty} \subset G(\R)$. At the finite places, we pick the open
compact subgroup $K_p=\Sp(n,\Z_p)$. Then $K = K_{\infty} \times
\prod_p K_p$ is the corresponding maximal compact subgroup of
$G(\A)$, the symplectic group over the ring of adeles of $\Q$. We
let $\omega = \omega_{n}$ be the Schr\"odinger model of the Weil
representation of $G_{\A}$ acting on $\calS(V^n(\A))$, the space of
Schwartz-Bruhat functions on $V^n(\A)$, associated to the standard
additive character of $\A/\Q$ (which on $\R$ is given by $t \mapsto
e(t)= e^{2\pi i t}$). Note that since $m$ is even we do not have to
deal with metaplectic coverings. We form the theta series associated
to $\varphi \in \calS(V^n(\A))$ by
\begin{align}
\theta(g,h,\varphi) = \sum_{\x \in V^n(\Q)}
(\omega(g)\varphi)(h^{-1}\x),
\end{align}
with $g \in G(\A)$ and $h \in H(\A)$. We assume $\varphi =
\varphi_{\infty} \otimes \varphi_f$ with $\varphi_{\infty} \in
\mathcal{S}(V^n(\R))$ and $\varphi_f \in \mathcal{S}(V^n(\A_f))$.

We now briefly review the Siegel-Weil formula, see e.g. \cite{Ku}.
We put
\begin{align}
I(g,\varphi) = \int_{H(\Q) \back H(\A)} \theta(g,h,\varphi) dh,
\end{align}
where $dh$ is the invariant measure on $H(\Q) \back H(\A)$
normalized to have total volume $1$. By Weil's convergence criterion
\cite{Weil}, $I(g,\varphi)$ is absolutely convergent if either $V$
is anisotropic or if
\begin{equation}\label{Weilcrit}
m -r > n+1.
\end{equation}
Here $r$ is the Witt index of $V$, i.e., the dimension of a maximal
isotropic subspace of $V$ over $\Q$.

We set $n(b)= \kzxz{1}{b}{0}{1}$ for $b$ a symmetric $n\times n$
matrix and $m(a) = \kzxz{a}{0}{0}{^ta^{-1}}$ for $ a \in \GL(n)$.
Then the Siegel parabolic is given by $P(\A) = N(\A)M(\A)$ with $N =
\{n(b);\; b \in \Mat_n,\, b ={^tb}\}$ and $M = \{m(a); \;a \in\GL(n)\}$.
Then using the Iwasawa decomposition $G(\A) = P(\A)K$ we define
\begin{equation}\label{Section}
\Phi(g,s)=
\left(\omega(g)\varphi\right)(0)\cdot\det|a(g)|_{\A}^{s-s_0},
\end{equation}
where
\begin{align}\label{s0}
s_0 = \frac{m}2- \frac{n+1}{2}.
\end{align}
Thus $\Phi$ defines a section in a certain induced parabolic induction
space (see \cite{Ku} (I.3.6)). Note that $\Phi$ is determined by its
values on $K$. Since $\Phi$ comes from $\varphi \in
\mathcal{S}(V^n(\A))$, we also see that $\Phi$ is a standard section,
i.e., its restriction to $K$ does not depend on $s$, and we write
$\Phi(k)= \Phi(k,s)$ for $k \in K$. Furthermore, $\Phi$ factors as $
\Phi = \Phi_{\infty} \otimes \Phi_f$.

%Furthermore, note that $\Phi(k)$ is left $\SO(n)$-invariant, which
%follows
%from the explicit formulas of the Weil representation. This implies that $\Phi(k)$ only depends on
%\[
%\det(k) \qquad \text{and} \qquad {^t k}k = \begin{pmatrix} a^2+c^2 & ab+cd \\ ab+cd & b^2 + d^2 \end{pmatrix},
%\]
%with $k = \left( \begin{smallmatrix} a&b \\c&d \end{smallmatrix}\right) \in \Uni(2)$.
We then define the Eisenstein series associated to $\Phi$ by
\begin{align}
E(g,s,\Phi) = \sum_{ \g \in P(\Q) \back G(\Q)} \Phi(\g g,s),
\end{align}
which for $\re(s) > \rho_n:=(n+1)/2$ converges absolutely and has a
meromorphic continuation to the whole complex plane. The extension
of Weil's work \cite{Weil} by Kudla and Rallis in the
convergent range is:

\begin{theorem}
\label{Siegel-Weil}
(\cite{KR1}, \cite{KR2}.)
Assume Weil's convergence criterion holds.
\begin{itemize}
\item[(i)]
Then $E(g,s,\Phi)$ is holomorphic at $s=s_0$.
\item[(ii)]
We have
\[
I(g,\varphi) = c_0 E(g,s_0,\Phi),
\]
where $c_0 =1$ if $m>n+1$ and $c_0 =2$ if $m \leq n+1$.

\end{itemize}

\end{theorem}

We translate the adelic Eisenstein series into more classical
language, see \cite{Ku} section IV.2. We let $K_f(N) \subset \prod_p
K_p$ be a subgroup of finite index of level $N$, i.e.,
\[
 \G := G(\Q) \cap (G(\R) K_f(N))
\]
contains the principal congruence subgroup $\G(N) \subset \Sp(n,\Z)$. We
assume that $\Phi_f$ is $K_f(N)$-invariant. Furthermore, if
$\varphi_f$ corresponds to the characteristic function of a coset of
an even lattice $L$ of level $N$ in $V$, then we have
\[
\Phi_f(\g) = \prod_{p |N} \Phi_p(\g)
\]
for $\g \in \G$. Via $G(\A)=G(\Q)G(\R)K_f(N)$ we see that the
Eisenstein series $E(g,s,\Phi)$ is determined by its restriction to
$G(\R)$. We assume that the restriction of $\Phi(g,s)$ to
$K_{\infty}$ is given by
\begin{align}\label{eq:phikappa}
\Phi_{\infty}^{\kappa}(k,s) := \det(\mathbf{k})^{\kappa}.
\end{align}
We denote the unique section at the Archimedian prime with this
property by $\Phi^{\kappa}_{\infty}$. Let $g_{\tau}= n(u)m(a)$ with
$^taa=v$ be an element moving the base point $i1_n$ of the
Siegel upper half plane $\h_n$ to $\tau = u+iv$. Then we obtain a
classical Eisenstein series of weight $\kappa$ (and level $N$):
\begin{align*}
E(g_{\tau},s,\Phi) & = \sum_{\g \in (P(\Q) \cap \G)\back \G}
\Phi^{\kappa}_{\infty}(\g g_{\tau}) \Phi_f(\g) \\
&= \det(v)^{\kappa/2} \sum_{\g \in (P(\Q) \cap \G)\back \G}
\left(\frac{\det(v)}{|\det(c\tau+d)|^2}\right)^{(s+\rho_n-\kappa)/2}
\det(c\tau+d)^{-\kappa} \; \Phi_f(\g),
\end{align*}
with $\g = \kzxz{a}{b}{c}{d}$. In particular, if $N=1$ then
\begin{align}
\label{classicaleis}
E(g_{\tau},s,\Phi) & =
\det(v)^{\kappa/2}  E^{(n)}_{\kappa}\big(\tau,(s+\rho_n-\kappa)/2\big),
\end{align}
where
\begin{align}
\label{classicaleis2}
E^{(n)}_{\kappa}(\tau,s)= \sum_{\g \in \G_{\infty} \back \Sp(n,\Z)}
\big(\det\Im(\gamma \tau)\big)^s
\det(c\tau+d)^{-\kappa}
\end{align}
is the classical Siegel Eisenstein series for $\Sp(n,\Z)$ of weight $\kappa$.

For later use, we introduce an embedding $\iota_0$ of $\Sp(n) \times
\Sp(n)$ into $\Sp(2n)$ by
\begin{equation}\label{embedding}
\zxz{a}{b}{c}{d} \times \zxz{a'}{b'}{c'}{d'} \mapsto
\begin{pmatrix}
a &  & b & \\
  &a'&   &b'\\
c &  & d &  \\
  &c'&   &d'
\end{pmatrix}.
\end{equation}

\section{Special Schwartz forms}\label{Schwartz}

We change the setting in this section and consider the real place
only. We assume that $V$ is now a real quadratic space of signature
$(p,q)$ of dimension $m$. Since it does not make any extra work we do
allow $m$ odd in this section. We set $H= \Orth(V)$. We pick an
oriented orthogonal basis $\{v_i\}$ of $V$ such that
$(v_{\alpha},v_{\alpha}) = 1$ for $\alpha = 1,\dots,p$ and
$(v_{\mu},v_{\mu}) = -1$ for $\mu = p+1,\dots,m$, and we denote the
corresponding coordinate functions by $x_{\alpha}$ and $x_{\mu}$. We
let $z_0$ be the $q$-dimensional subspace $\Span\{v_{p+1},\dots,v_m\}$
with the induced orientation. We let $K^{H}$ be the maximal compact
subgroup of $H$ stabilizing $z_0$. We realize the symmetric space
$D=H/K^H$ associated to $V$ as the Grassmannian of oriented negative
$q$-planes in $V$. Thus $D$ has two components
\[
D =D_+ \amalg D_-,
\]
where
\begin{equation*}
D_+ = \{ z \subset  V;\quad\text{$\dim z =q$, $(\,,\,)|_z < 0$,
$z$ has the same orientation as $z_0$} \}.
\end{equation*}
Thus $D_+ \simeq H_0/K^{H_0}$ with $H_0 = \SO_0(V)$, the connected
component of the orthogonal group, and $K^{H_0}$ the maximal compact
subgroup of $H_0$ stabilizing $z_0$. We associate to $z \in D$ the
standard majorant $(\,,\,)_z$ given by
\begin{equation*}
(x,x)_z = (x_{z^{\perp}},x_{z^{\perp}}) - (x_z,x_z),
\end{equation*}
where $x = x_z + x_{z^{\perp}} \in V$ is given by the orthogonal
decomposition $V = z \oplus z^{\perp}$. We write $(\,,\,)_0 =
(\,,\,)_{z_0}$.

Let $\mathfrak{h}$ be the Lie algebra of $H$  and let
$\mathfrak{h} = \mathfrak{p} \oplus \mathfrak{k}$ with
$\mathfrak{k}^H = \operatorname{Lie}(K^{H})$ be the associated
Cartan decomposition. Then $\mathfrak{p} \simeq \mathfrak{h}/
\mathfrak{k}^H$ is isomorphic to the tangent space at the base point
$z_0$ of $D$. With respect to the above basis of $V$ we have
\begin{equation}\label{mathfrakp}
\mathfrak{p} \simeq \left\{
\begin{pmatrix}
0 & X \\
{^tX} &0
\end{pmatrix}
;\,\, X \in \Mat_{p,q}(\R) \right\}.
\end{equation}
We let $X_{\alpha\mu}$ ($ 1 \leq \alpha \leq p$, $p+1 \leq \mu \leq
p+q$) denote the element of $\mathfrak{p}$ which interchanges
$v_{\alpha}$ and $v_{\mu}$ and annihilates all the other basis
elements of $V$. We write $\omega_{\alpha\mu}$ for the element of
the dual basis corresponding to $X_{\alpha\mu}$.

We let $\omega=\omega_n$ be the Weil representation of the metaplectic
cover $\Mp(n,\R)$ of $\Sp(n,\R)$ acting on the Schwartz functions
$\mathcal{S}(V^n)$. We let
$K=\widetilde{\Uni}(n)$ be the maximal compact subgroup of
$\Mp(n,\R)$ given by the inverse image of the standard maximal compact subgroup $\Uni(n)$ in $\Sp(n,\R)$.
%
%$K_0\simeq \Uni(n)$ be the maximal compact
%in $\Sp(n,\R)$ and let $K=\widetilde{K}_0$ be the maximal compact
%subgroup in $\Mp(n,\R)$ given by the full inverse image of $K_0$ in
%$\Mp(n,\R)$.
Recall that $K$ admits a character $\det^{1/2}$ whose
square descends to the determinant character of $\Uni(n)$. We also
write $\omega$ for the associated Lie algebra action on the space of
$K$-finite vectors in $\mathcal{S}(V^n)$. It is given by the so-called
polynomial Fock space $S(V^n) \subset \mathcal{S}(V^n)$. It consists
of those Schwartz functions on $V^n$ of the form $p(\x)\varphi_0(\x)$,
where $p(\x)$ is a polynomial function on $V^n$. Here $\varphi_0(\x)$
is the standard Gaussian on $V^n$. More precisely, for $\x
=(x_1,\dots,x_n) \in V^n$ and $z\in D$, we let
\[
\varphi_0(\x,z) = \exp\left(-\pi \sum_{i=1}^n (x_i,x_i)_z\right),
\]
and set $\varphi_0(\x) = \varphi_0(\x,z_0)$. We view
\begin{equation}
\varphi_0 \in  [\mathcal{S}(V^n) \otimes C^{\infty}(D)]^{H} \simeq
[\mathcal{S}(V^n) \otimes {\bigwedge}
^0(\mathfrak{p^{\ast}})]^{K^{H}},
\end{equation}
where the isomorphism is given by evaluation at the base point $z_0$
of $D$. In the following we will identify corresponding objects
under this isomorphism.

Kudla and Millson (see \cite{KMI}) constructed (in much greater
generality) Schwartz forms $\varphi_{KM}$ on $V$ taking values in
$\mathcal{A}^q(D)$, the differential $q$-forms on $D$. More
precisely,
\begin{equation*}
\varphi_{KM} \in [ \mathcal{S}(V) \otimes \mathcal{A}^q(D)]^{H}
\simeq
 [ \mathcal{S}(V) \otimes {\bigwedge} ^q(\mathfrak{p^{\ast}})]^{K^{H}},
\end{equation*}
where the isomorphism is again given by evaluation at the base point
of $D$. The Schwartz form $\varphi_{KM}$ is given by
\begin{gather*}
\varphi_{KM}= \frac{1}{2^{q/2}} \prod_{\mu = p+1}^{p+q} \left[
\sum_{\alpha =1}^{p} \left(  x_{\alpha} - \frac{1}{2\pi}
\frac{\partial}{\partial x_{\alpha}}  \right)   \otimes
A_{\alpha\mu}  \right] \varphi_0.
\end{gather*}
Here $A_{\alpha\mu}$ denotes the left multiplication by
$\omega_{\alpha\mu}$. More generally, we consider the Schwartz forms
\[
\varphi_{q,\ell} \in [\mathcal{S}(V) \otimes {\bigwedge}
^q(\mathfrak{p^{\ast}}) \otimes \Sym^{\ell}(V)]^{K^{H}}
\]
with values in the $\ell$-th symmetric powers of $V$ introduced by
Funke and Millson \cite{FM}. Viewing the Kudla-Millson form
$\varphi_{KM}$ as the form $\varphi_{q,0} \in [\mathcal{S}(V)
\otimes {\bigwedge} ^q(\mathfrak{p^{\ast}}) \otimes
\Sym^0(V)]^{K_{H}}$, the forms $\varphi_{q,\ell}$ are given by
\begin{align*}
\varphi_{q,\ell}& = \left[\frac12 \sum_{\alpha=1}^p \left(
x_{\alpha} - \frac{1}{2\pi} \frac{\partial}{\partial x_{\alpha}}
\right) \otimes 1 \otimes A_{v_{\alpha}} \right]^{\ell} \varphi_{KM}
\\ &= \frac{1}{2^{\ell}}\sum_{\alpha_1,\dots,\alpha_{\ell}=1}^{p} \left[
\prod_{i=1}^{\ell} \left( x_{\alpha_i} - \frac{1}{2\pi}
\frac{\partial}{\partial x_{\alpha_i}} \right) \otimes 1 \otimes
\prod_{i=1}^{\ell} A_{v_{\alpha_i}} \right] \varphi_{KM}.
\end{align*}
Here $A_v$ denotes the multiplication with the vector $v$ in the
symmetric algebra of $V$. Note that $\Sym^{\ell}(V)$ is not an
irreducible representation of $H$, and we denote by
$\varphi_{q,[\ell]}$ the projection of $\varphi_{q,\ell}$ onto
$\mathcal{H}^{\ell}(V)$, the harmonic $\ell$-tensors in $V$. It
consists of those symmetric $\ell$-tensors which are annihilated by
the signature $(p,q)$-Laplacian $\Delta= \sum_{\alpha =1}^p
\tfrac{\partial^2}{\partial v_{\alpha}^2}- \sum_{\mu=p+1}^m
\tfrac{\partial^2}{\partial v_{\mu}^2}$. Here we view $v_{\alpha}$
and $v_{\mu}$ as independent variables. It can be also characterized
as the space of symmetric $\ell$-tensors in $V$ which are orthogonal
with respect to the induced inner product on $\Sym^{\ell}(V)$ to
vectors of the form $r^2 w$. Here $w \in \Sym^{\ell-2}(V)$ and $r^2$
denotes the multiplication with $\sum_{\alpha=1}^p v_{\alpha}^2 -
\sum_{\mu = p+1}^m v^2_{\mu}$. Recall that we have
$\Sym^{\ell}(V)=\mathcal{H}^{\ell}(V) \oplus r^2 \Sym^{\ell-2}(V)$
as representations of $H$.

The Schwartz form $\varphi_{q,\ell}$ (and also $\varphi_{q,[\ell]}$)
is an eigenfunction of weight $m/2+\ell$ under the action of $k \in
K$, see \cite{KMI,FM}, i.e.,
\begin{equation}\label{trafo0}
\omega(k) \varphi_{q,\ell}= \det(\mathbf{k})^{m/2+\ell}
\varphi_{q,\ell}.
\end{equation}
Here $\mathbf{k}$ is the element in $\widetilde{\Uni}(1)$ corresponding
to $k \in \widetilde{\SO}(2)\subset \Mp(1,\R)$.
Moreover, $\varphi_{q,\ell}(x)$ is a closed differential form on $D$.

We normalize the inner product on $\Sym^{\ell}(V)$ inductively by
setting
\[
(w_1\cdots w_{\ell},w'_1\cdots w'_{\ell}) = \frac1{\ell}
\sum_{j=1}^{\ell} (w_1,w'_j)(w_2\cdots w_{\ell}, w'_1 \cdots
\widehat{w'_j}\cdots w'_{\ell}).
\]
With this normalization we easily see that for the
restriction of $(\,,\,)$ to the positive definite subspace
$\Span\{v_{\alpha};\; 1 \leq \alpha \leq p\}$ of $V$ we have
\[
\sum_{\substack{\alpha_1,\dots,\alpha_{\ell}=1
\\\beta_1,\dots,\beta_{\ell}=1} }^p \left(\prod_{i=1}^{\ell}
v_{\alpha_i}, \prod_{i=1}^{\ell} v_{\beta_i}\right) = p^{\ell}.
\]
We let $\widetilde{\Sym}^{\ell}(V)$ be the local system on $D$
associated to $\Sym^{\ell}(V)$. Then for the wedge product, we have
$\wedge: \mathcal{A}^{r}(D,\widetilde{\Sym}^{\ell}(V))
\times\mathcal{A}^{s}(D,\widetilde{\Sym}^{\ell}(V)) \to
\mathcal{A}^{r+s}(D)$ by taking the inner product on the fibers
$\Sym^{\ell}(V)$. We are ultimately more interested in the form
$\varphi_{q,[\ell]}$, but calculations with $\varphi_{q,\ell}$ are
more convenient. In this context the following lemma will be
important later.

\begin{lemma}\label{Brauer}
Let $\eta \in \mathcal{A}^{(p-1)q}(D,\widetilde{\Sym}^{\ell-2}(V))$.
Then
\[
\varphi_{q,\ell} \wedge r^2 \eta = -\frac{1}{2\pi} (\omega(R)
\varphi_{q,\ell-2}) \wedge \eta.
\]
Here $R = \tfrac12 \left( \begin{smallmatrix} 1&i \\ i &-1
\end{smallmatrix}\right) \in \mathfrak{sl}(2,\C)$ is the standard
$\SL(2)$-raising operator.
\end{lemma}

\begin{proof}
By the adjointness of $\tfrac1{\ell(\ell-1)}\Delta$ and $r^2$ with
respect to the inner product in $\Sym^{\bullet}(V)$, we have
$\varphi_{q,\ell} \wedge r^2 \eta = \tfrac1{\ell(\ell-1)} (\Delta
\varphi_{q,\ell}) \wedge \eta$. Note that $\Delta$ operates on the
coefficient part of $\varphi_{q,\ell}$. Then switching to the Fock
model of the Weil representation, see the proof of
Lemma~\ref{phiformel}, and using \eqref{raisingaction} one easily
sees $\Delta \varphi_{q,\ell} = - \tfrac{\ell(\ell-1)}{2\pi}
\omega(R) \varphi_{q,\ell-2}$. We leave the details to the reader.
\end{proof}

We let $\ast$ denote the Hodge $\ast$-operator on $D$. Then
$\varphi_{q,\ell}(x_1) \wedge \ast \varphi_{q,\ell}(x_2)$ with
$\x=(x_1,x_2) \in V^2$, being a top degree differential form, gives
rise to a scalar-valued Schwartz function $\phi_{q,\ell}$ on $V^2$
defined by
\begin{equation}\label{phiKMdef}
\phi_{q,\ell}(\x,z) \mu = \varphi_{q,\ell}(x_1,z) \wedge  \ast
\varphi_{q,\ell}(x_2,z).
\end{equation}
Here $\mu$ is the volume form on $D$ induced by the Riemannian
metric coming from the Killing form on $\mathfrak{g}$. For
convenience we scale the metric such that the restriction of $\mu$ to
the base point $z_0$ is given by
\begin{align}
\label{metric-normalization}
\mu = \omega_{1, p+1} \wedge \cdots \wedge \omega_{1,p+q} \wedge
\omega_{2,p+1} \wedge \cdots \wedge \omega_{p,p+q}.
\end{align}
Note that
\[
\phi_{q,\ell} \in  [\mathcal{S}(V^2) \otimes C^{\infty}(D)]^{H_0}
\simeq [\mathcal{S}(V^2) \otimes {\bigwedge}
^0(\mathfrak{p^{\ast}})]^{K^{H_0}}.
\]

\begin{lemma}
We have
%For a multi-index $\underline{\alpha} =(\alpha_1,\dots, \alpha_q) \in \{1,\dots,p\}^q$,
\begin{align*}
\phi_{q,\ell}(\x) &=  \frac{p^{\ell}}{2^{q+2\ell}}
\sum_{\alpha_1,\dots,\alpha_{q+\ell}=1}^{p}
\prod_{i=1}^{q+\ell}\left( x_{\alpha_i 1} - \frac{1}{2\pi}
\frac{\partial}{\partial x_{\alpha_i 1}} \right) \left(  x_{\alpha_i
2} - \frac{1}{2\pi}
\frac{\partial}{\partial x_{\alpha_i 2}}  \right) \varphi_0(\x) \\
 &= \frac{p^{\ell}}{2^{q+2\ell}} \left( \sum_{\alpha=1}^p
 \left( x_{\alpha 1} - \frac{1}{2\pi} \frac{\partial}{\partial x_{\alpha 1}}  \right) \left(  x_{\alpha 2} - \frac{1}{2\pi}
\frac{\partial}{\partial x_{\alpha 2}}  \right)
\right)^{q+\ell}\varphi_0(\x).
\end{align*}
\end{lemma}

\begin{example}
%For signature $(p,1)$, we have
%\[
%\phi_{q,0}(\x) = 2 \sum_{\alpha=1}^p x_{\alpha 1} x_{\alpha 2}
%\varphi_0(\x),
%\]
For signature $(p,2)$, we have
\[
\phi_{q,0}(\x) = \sum_{\alpha=1}^p (x_{\alpha 1}^2 - \frac{1}{4\pi}) (x_{\alpha 2}^2 - \frac{1}{4\pi}) \varphi_0(\x) + 4\sum_{\substack{\alpha,\beta=1\\
\alpha \ne \beta}}^p x_{\alpha 1} x_{\beta 1} x_{\alpha 2} x_{\beta
2} \varphi_0(\x).
\]
\end{example}

Note that \eqref{trafo0} immediately implies:

\begin{lemma}
For $k_1,k_2 \in \widetilde{\SO}(2) \subset \Mp(1,\R)$, we have
\[
\omega(\iota_0(k_1,k_2)) \phi_{q,\ell} =
\det(\mathbf{k_1}\mathbf{k_2})^{m/2+\ell} \phi_{q,\ell}.
\]
\end{lemma}

The action of the full maximal compact $K \subset \Mp(2,\R)$ on
$\phi_{q,\ell}$ via the Weil representation is more complicated, as
we now explain.
We let
\begin{equation}
\mathfrak{g} = \mathfrak{k} \oplus \mathfrak{p}_+\oplus
\mathfrak{p}_-
\end{equation}
be a Harish-Chandra decomposition of $\mathfrak{g} =
\mathfrak{sp}(2,\C)$, where $\mathfrak{k}= \Lie(K)_{\C}$,
\begin{equation}
\mathfrak{p}_+ = \left\{ p_+(X) = \frac12 \begin{pmatrix} X&  iX \\
iX & -X\end{pmatrix}; \; X \in \Mat_2(\C), {^tX}=X \right\},
\end{equation}
and $\mathfrak{p}_- = \overline{\mathfrak{p_+}}$. Note that
$\mathfrak{p}_+$ is the holomorphic tangent space of $\h_2$ at the
base point $i1_2$ and is spanned by the raising operators
\begin{gather}
R_{1}= R_{11}= p_+ \left(\begin{smallmatrix} 1&0 \\0 &
0\end{smallmatrix}\right), \qquad R_2 = R_{22} =p_+ \left(
\begin{smallmatrix} 0& 0\\0 &1 \end{smallmatrix}\right), \\
R_{12} = \frac12 p_+ \left(  \begin{smallmatrix}0 &1 \\1 &0
\end{smallmatrix} \right).
\end{gather}
Note that $R_1 = \iota_0(R,0)$ and $R_2 = \iota_0(0,R)$ are the
images of the $\SL_2$-raising operator $R$ in $\mathfrak{sp}(2,\C)$
under the two standard embeddings of $\mathfrak{sl}(2)$ into
$\mathfrak{sp}(2)$.

Recall that the adjoint action of $K$ on $\mathfrak{p}_+$ is
isomorphic to the standard action of $K$ on $\Sym^2(\C^2)$.
Explicitly, the intertwiner is given by $R_{rs} \mapsto e_r e_s$,
where $e_1, e_2$ denotes the standard basis of $\C^2$. We obtain an
isomorphism of $K$-modules
\begin{equation}
\Sym^{\bullet}\Sym^2\C^2 = \bigoplus_{j=0}^{\infty}
\Sym^{j}\Sym^2\C^2 \simeq U(\mathfrak{p}_+)
\end{equation}
of the symmetric algebra on $\Sym^2\C^2$ with the universal
enveloping algebra of $\mathfrak{p}_+$.

\begin{lemma}\label{phiformel}
We have
\[
\phi_{q,\ell} = \frac{p^{\ell}(-1)^{q+\ell}}{2^{\ell}\pi^{q+\ell}}
\omega(R_{12})^{q+\ell} \varphi_0.
\]
\end{lemma}

\begin{proof}
We indicate a quick proof using the Fock model of the Weil
representation. For more details for what follows, see the appendix
of \cite{FM}. There is an intertwining map $\iota: S(V^n) \to
\mathcal{P}(\C^{n(p+q)})$ from the  polynomial Fock space to the
infinitesimal Fock model of the Weil representation acting on the
space of complex polynomials $\mathcal{P}(\C^{n(p+q)})$ in $n(p+q)$
variables such that $\iota(\varphi_0) = 1$. We denote the variables
in $\mathcal{P}(\C^{n(p+q)})$ by $z_{\alpha i}$ ($1 \leq \alpha \leq
p$) and $z_{\mu i}$ ($p+1 \leq \mu \leq p+q$) with $i=1,\dots ,n$.
Moreover, the intertwining map $\iota$ satisfies
\[
\iota \left(  x_{\alpha i} - \frac{1}{2\pi} \frac{\partial}{\partial
x_{\alpha i}}  \right)  \iota^{-1} =
 \frac{1}{2\pi i} z_{\alpha i}.
\]
Hence in the Fock model, we have
\[
\phi_{q,\ell} = \frac{p^{\ell}}{2^{q+2\ell}} \left(\frac{1}{2\pi i}
\right)^{2(q+\ell)} \left[ \sum_{\alpha=1}^p z_{\alpha 1 } z_{\alpha
2} \right]^{q+\ell}.
\]
On the other hand, for the action of the raising operators, we find
\begin{equation}\label{raisingaction}
\omega(R_{rs}) =   \frac1{8\pi} \sum_{\alpha=1}^p z_{\alpha r }
z_{\alpha s} -2 \pi \sum_{\mu=p+1}^m \frac{\partial^2}{\partial
z_{\mu r}\partial z_{\mu s}}.
\end{equation}
In the Fock model, we therefore have $\omega(R_{12})^{q+\ell}
\varphi_0 = \left[\frac1{8\pi} \sum_{\alpha=1}^p z_{\alpha 1 }
z_{\alpha 2} \right]^{q+\ell}$, and the lemma follows.
%For completeness, $\mathfrak{k}$ acts linearly on the polynomials.
%More precisely, we let $w_{rs} \in \mathfrak{k}$ be the endomorphism
%of $\C^2$ mapping $e_r$ to $e_s$ and annihilating the other basis
%element. Then the Fock model action of $w_{rs}$ is given by
%\begin{equation}
%\omega(w_{rs})=\left[ \sum_{\alpha=1}^p z_{\alpha s} \frac
%{\partial}{\partial z_{\alpha r}} - \sum_{\mu = p+1}^{p+q} z_{\mu r}
%\frac {\partial}{\partial z_{\mu s}} \right] + \frac12(p-q)
%\delta_{rs}.
%\end{equation}
\end{proof}

We obtain:

\begin{proposition}\label{trafo}
For $k \in K \simeq \widetilde{\Uni}(2)$, we have
\[
\omega(k) \phi_{q,\ell}
=\frac{p^{\ell}(-1)^{q+\ell}}{2^{\ell}\pi^{q+\ell}}
\det(\mathbf{k})^{(p-q)/2} \left( \Ad(k) R_{12} \right)^{q+\ell}
\varphi_0.
\]
\end{proposition}

\begin{proof}
This follows immediately from Lemma~\ref{phiformel} and the fact
that the Gaussian $\varphi_0$ has weight $(p-q)/2$.
\end{proof}

\begin{remark}
The Kudla-Millson forms $\varphi_{KM}$ cannot be expressed in terms
of elements in $\mathfrak{p}_+$.
\end{remark}

Proposition~\ref{trafo} reduces the $K$-action on $\phi_{q,\ell}$ to
the representation theory of the group $\Uni(2)(\C)=\GL_2(\C)$ on
$\Sym^{\bullet}\Sym^2\C^2$, which is given as follows.

\begin{lemma}
\label{lem:repsplit}
The $\GL_2(\C)$-representation $\Sym^{j}\Sym^2\C^2$
decomposes as
\[
\Sym^{j}\Sym^2\C^2 \simeq \bigoplus_{i=0}^{[j/2]} \Sym^{2j-4i}\C^2
\otimes \det{^{2i}}
\]
into its irreducible constituents. The summand for $i = [j/2]$ is
given by
\begin{equation}\label{rep1}
 \Sym^{2j-4[j/{2}]} \C^2 \otimes \det{^{2[j/{2}]}} =
\begin{cases} \det {{^{j}}} & \text{if $j$ is even}, \\
              \Sym^2 \C^2 \otimes \det{^{j-1}}  & \text{if $j$ is odd},
\end{cases}
\end{equation}%\label{vector}
and is generated by the vector
\begin{equation}\label{rep2}
\alpha_{j} = \sum_{i=0}^{[j/2]} \left( \begin{smallmatrix} [j/2]
\\ i  \end{smallmatrix} \right)(-1)^i (e_1^{2})^i (e_2^{2})^i
(e_1e_2)^{j-2i} = \begin{cases} \left[ (e_1e_2)^2 - e_1^2e_2^2
\right]^{j/2}& \text{if $j$ is even}, \\
(e_1e_2)\left[ (e_1e_2)^2 - e_1^2e_2^2 \right]^{[j/2]}& \text{if $j$
is odd}.
\end{cases}
\end{equation}

\end{lemma}

\begin{proof}
For the first statement, see e.g. \cite{FH}, p.81/82. For
\eqref{rep2}, note that in
\[\Sym^2\Sym^2\C^2 = \Sym^4 \C \oplus
\det{^2},\]
the vector
\[
\alpha_2 =  (e_1e_2)^2 - e_1^2e_2^2
\]
generates the one-dimensional sub-representation. Then, for $j$
even, $\alpha_j$ is given by the image of $\left( \alpha_2
\right)^{j/2} \in \Sym^{j/2}\Sym^2\Sym^2 \C^2$ under the projection
onto $\Sym^{j}\Sym^2 \C^2$. The argument for $j$ odd is analogous.
\end{proof}

By slight abuse of notation, we also write $\alpha_j$ for the
corresponding element in $\Uni(\mathfrak{p}_+)$ and define another
Schwartz function $\xi=\xi_{q,\ell} \in \calS(V^2)$ by
\begin{equation}\label{xi:def}
\xi =
\xi_{q,\ell}=\frac{p^{\ell}(-1)^{q+\ell}}{2^{\ell}\pi^{q+\ell}}
\omega(\alpha_j) \varphi_0.
\end{equation}

\begin{proposition}\label{decomp}
For the Schwartz function $\phi_{q,\ell}$, there exists a $\psi \in
\mathcal{S}(V^2)$ such that
\begin{align}
\phi_{q,\ell} = \xi_{q,\ell}+ \omega(R_1)\omega(R_2) \psi.
\end{align}
\end{proposition}

\begin{proof}
We have
\[
(e_1e_2)^{q+\ell} - \alpha_{q+\ell} = e_1^2e_2^2
\sum_{i=1}^{[(q+\ell)/2]} \left(
\begin{smallmatrix}  [(q+\ell)/2] \\ i  \end{smallmatrix} \right)(-1)^i
(e_1^{2})^{i-1} (e_2^{2})^{i-1} (e_1e_2)^{q+\ell-2i}.
\]
Using the intertwiner with $\Uni(\mathfrak{p}_+)$, we recall that
$e_i^2$ corresponds to $R_{i}$. Thus $\psi$ is given by
\[
\psi =
\frac{p^{\ell}(-1)^{q+\ell}}{2^{\ell}\pi^{q+\ell}}\sum_{i=1}^{[(q+\ell)/2]}
\left(
\begin{smallmatrix}  [(q+\ell)/2] \\ i  \end{smallmatrix} \right)(-1)^i
\omega \left(R_1^{i-1} R_2^{i-1} R_{12}^{q+\ell-2i} \right)
\varphi_0.
\]
\end{proof}

One easily sees using \eqref{raisingaction}:

\begin{lemma}\label{xi:vanishing}
The Schwartz function $\xi$ vanishes identically if and only if
$p=1$ and $q+\ell>1$.
\end{lemma}

\begin{example}
%For $q$ even, $p>1$, and $\ell=0$, we actually have
%\[
%\xi e_q^p= C \varphi_{KM} \wedge \varphi_{KM} \wedge e_q^{p-2}
%\]
%for a nonzero constant $C$ and where $e_q$ denotes the Euler form on
%$D$.
For $q=2$, $p>1$, and $\ell=0$, we have
\[
\phi_{2,0} \cdot \Omega^{p} = C  \varphi_{KM} \wedge \varphi_{KM} \wedge
\Omega^{p-2} + C' \omega(R_1)\omega(R_2) \varphi_0 \cdot \Omega^{p}
\]
for some nonzero constants $C$ and $C'$. Here $\Omega$ denotes the
K\"ahler form on the Hermitian domain $D$. But we will not need this.
\end{example}

In view of Lemma \ref{lem:repsplit} and Proposition \ref{trafo}, we see for $q+\ell$ even that
\begin{equation}\label{trafoxi}
\omega(k) \xi = \det(\mathbf{k})^{m/2+\ell} \xi
\end{equation}
for $k\in K$. We let $\Xi(g,s)$ be the section in
the induced representation corresponding to the Schwartz function $\xi$ via \eqref{Section}.

\begin{proposition}\label{Xi:Indent}
Let $q +\ell$ be even. Then $\Xi$ is the standard section
\eqref{eq:phikappa} at the infinite place of weight $m/2+\ell$. More
precisely,
\begin{equation}\label{Xi:eq}
\Xi(s)= C(s)\Phi_{\infty}^{m/2+\ell}(s)
\end{equation}
for a certain (explicit) polynomial $C(s)$. Moreover,
\[
C(s_0) \ne 0
\]
with $s_0 = (m-3)/2$ as in \eqref{s0} for $p>1$, while $C(s) \equiv
0$ for $p=1$.
\end{proposition}

\begin{proof}
The identity \eqref{Xi:eq} follows from \eqref{trafoxi} and the
uniqueness of $\Phi_{\infty}^{m/2+\ell}$. The precise statement
follows from considerations in \cite{KR}. The element
$\alpha_{q+\ell}$ is trivially a highest weight vector of weight $\mu=(q+\ell,q+\ell)$ of $\GL_2(\C)$. Therefore we can take $\alpha_{q+\ell}$  equal to the element $u_{\mu}^0 \in
\Uni(\mathfrak{p}_+)$ (or $u_{\mu} \in \Uni(\mathfrak{g})$) in the
notation of \cite{KR}, p.31/32. Then by Corollary~1.4 of \cite{KR}, we have
$\Xi(s)= u_{\mu} \Phi_{\infty}^{(p-q)/2}(s) = cP_{\mu}^{(p-q)/2}(s)
\Phi^{m/2+\ell}(s)$, for a certain polynomial $P_{\mu}^{(p-q)/2}$
and a nonzero constant $c$. One easily sees $P_{\mu}^{(p-q)/2}(s_0)
\ne 0$ for $p>1$. See also \cite{KR}, p. 38. For $p=1$, $\Xi$
vanishes identically, since already $\xi=0$ by
Lemma~\ref{xi:vanishing}.
\end{proof}

\begin{remark}
For $q+\ell$ odd, we see in the same way
\[
\Xi(s) = C(s) R_{12} \Phi_{\infty}^{m/2+\ell-1}(s)
\]
for a certain polynomial $C(s)$. Note that $\alpha_{q+\ell}$ is
\emph{not} a highest weight vector for $\Sym^2 \C^2 \otimes
\det{^{q+\ell-1}}$ (which has weight $(q+\ell+1,q+\ell-1)$).
\end{remark}

\section{The $L^2$-norm of the theta lift}
\label{sect:5}

We now return to the global situation and retain the notation of
Section 2. Let $V$ be a non-degenerate quadratic space over $\Q$ of
signature $(p,q)$ and even dimension $m=p+q$.
We let $L\subset V$ be an even lattice
and write $L^{\#}$ for the dual lattice. For each prime $p$,
we let $L_p =L \otimes \Z_p$ and let $K^{H}_p$ be the subgroup of
$\Orth(L_p)$ given by the kernel of $\Orth(L_p) \rightarrow
\Orth(L_p^{\#}/L_p)$. Then $K_f^H = \prod_p K_p^H$ is an open
compact subgroup of $H(\A_f)$. We write $H(\R)_0 = \SO_0(V(\R))$,
and we let $K^H_{\infty}$ be a maximal compact subgroup of $H(\R)$.
Then $D = H(\R) / K^H_{\infty}$ is the symmetric domain of oriented
negative $q$-planes considered in the previous section.
%We let $H(\Q)_0 =
%H(\Q) \cap H(\R)_0H(\A_f)$, and
By strong approximation we write
\begin{align}
\label{strong-approx}
H(\A) = \coprod_j H(\Q) H(\R)_0 h_j K_f^H
\end{align}
with $h_j \in H(\A_f)$. Then we put
\begin{align}
X = X_{K_f^H}=H(\Q) \back (D \times H(\A_f))/K_f^H
\end{align}
such that
\begin{align}
X \simeq \coprod_j X_j
\end{align}
with $X_j = \G_j \back D_+$, where $\G_j = H(\Q) \cap (H(\R)_0
h_jK_f^Hh_j^{-1}).$ We let $\varphi_f \in
\mathcal{S}(V(\A_f))^{K_f^H}$ be a $K_f^H$-invariant Schwartz
function on the finite adeles. Then $\varphi_f$ corresponds to a
linear combination of characteristic functions on the discriminant
group $L^{\#}/L$. Since $\varphi_{q,\ell}$ is an eigenfunction of
weight
\[
\kappa =m/2+\ell
\]
under the action of $\Uni(1)$, we can
form the classical theta function on $\h$, the upper half space, by
setting
\begin{align*}
\theta(\tau,z,\varphi_{q,\ell}\otimes \varphi_f) &= v^{-\kappa/2}
\sum_{x \in V(\Q)}\varphi_f(x) \omega_{\infty}(g_{\tau})
\varphi_{q,\ell}(x,z) \\ &= v^{-\ell/2} \sum_{x \in V(\Q)}
\varphi_f(x) \varphi_{q,\ell}(\sqrt{v}x,z) e^{\pi i (x,x)u}.
\end{align*}
Here $\tau = u+iv \in \h$, and $g_{\tau} = \kzxz{1}{u}{0}{1}
\kzxz{\sqrt{v}}{0}{0}{\sqrt{v}^{-1}} \in G(\R) \subset G(\A)$ is the
standard element moving the base point $i \in \h$ to $\tau$. Then
$\theta(\tau,z,\varphi_{q,\ell}\otimes\varphi_f)$ transforms like a
non-holomorphic modular form of weight $\kappa$ for the principal
congruence subgroup $\G(N)$ of $\SL_2(\Z)$ taking values in the
differential $q$-forms on $X$. Here $N$ is the level of $L$, i.e.,
the smallest positive integer such that $\tfrac12 N(x,x) \in \Z$ for
all $x \in L^{\#}$. In particular, if $L$ is unimodular,
$\theta(\tau,z,\varphi_{q,\ell}\otimes \varphi_f)$ is a form for the
full modular group $\SL_2(\Z)$.

%\begin{remark}
%For $q+\ell$ odd and $\beta=0$, the theta series
%$\theta(\tau,\varphi_{q,\ell})$ vanishes identically. To obtain
%non-zero examples in this case, one needs the coset condition.
%\end{remark}

We write $S_{\kappa}(\G(N))$ for the space of cusp forms of weight
$\kappa$ for $\G(N)$. We normalize the Petersson scalar product be putting
\begin{align}
\label{petnorm}
(f,g)=\frac{1}{[\Gamma(1):\Gamma(N)]}\int_{\Gamma(N)\bs \H} f(\tau) \overline{g(\tau)} v^{\kappa}\, d\mu(\tau)
\end{align}
for $f,g\in S_{\kappa}(\G(N))$.
Here $ d\mu(\tau) = \tfrac{du\,dv}{v^2}$ is the invariant measure on
$\h$.
For $f \in S_{\kappa}(\G(N))$, we consider the theta lift
\begin{align}
\Lambda(f) = \big(f,\theta(\tau,\varphi_{q,\ell}\otimes\varphi_f)\big) =
\int_{\G(N) \back \h} f(\tau)
\overline{\theta(\tau,\varphi_{q,\ell}\otimes\varphi_f)} v^{\kappa}\,
d\mu(\tau).
\end{align}
%
%We let $\mathcal{H}^q(X)$ be the space of harmonic $q$-forms on $X$.
%Then
%\begin{proposition}[\cite{KMCan}, Theorem 4.1.]
%Let $f \in S_{\kappa}(\G(N))$. Then
%\[
%\Lambda(f) \in \mathcal{H}^q(X),
%\]
%i.e, we have a map
It defines  a linear map
\begin{align}
\Lambda:S_{\kappa}(\G(N)) \longrightarrow
\mathcal{Z}^q(X,\widetilde{\Sym}^{\ell}(V))
\end{align}
into the $\widetilde{\Sym}^{\ell}(V)$-valued closed differential $q$-forms
on $X$.
%\end{proposition}

%\texttt{Der Lift sollte auch im allgemeinen mittels eines Einzeilers
%harmonisch sein, aber das bekomme ich gerade nicht so recht hin}

In order to show the injectivity of $\La$, we study its $L^2$-norm
given by
\begin{align}
\|\Lambda(f)\|_2^2 = \int_X\Lambda(f) \wedge \ast
\overline{\Lambda(f)}.
\end{align}
We will use the \emph{doubling method} to compute
$\|\Lambda(f)\|_2^2$, see \cite{Boe,Ga,PS-R,Ra}.

\begin{proposition}\label{doubling0}
Assume that $m>3+r$ so that Weil's convergence
criterion~\eqref{Weilcrit} in genus $2$ holds. Then $\Lambda(f)$ is
square integrable, and
\begin{equation}\label{doubling1}
\|\Lambda(f)\|_2^2 = \left(f(\tau_1) \otimes \overline{f(\tau_2)},
\tilde{I}(\tau_1,-\bar{\tau}_2,\phi_{q,\ell}\otimes \phi_f)\right),
\end{equation}
where $(\,,\,)$ denotes the Petersson scalar product on $\G(N)
\times \G(N)$ and
\begin{equation}\label{intdef1}
\tilde{I}(\tau_1,{\tau}_2,\phi_{q,\ell}\otimes\phi_f) = \int_X
\theta(\tau_1,{\tau}_2,z,\phi_{q,\ell}\otimes\phi_f) \mu
\end{equation}
is the integral over the locally symmetric space of the theta series
\begin{equation}\label{thetadef1}
\theta(\tau_1,{\tau}_2,z,\phi_{q,\ell}\otimes \phi_f) =
(v_1v_2)^{-\kappa/2} \sum_{\x \in V^2(\Q)} \phi_f(\x)
(\omega_{\infty}(\iota_0(g_{\tau_1},g_{\tau_2}))\phi_{q,\ell}(\x,z),
\end{equation}
(which by \eqref{trafo0} defines a modular form of weight $\kappa$
on $\G(N) \times \G(N)$). Here $\phi_f = \varphi_f \otimes \varphi_f
\in \mathcal{S}(V^2(\A_f))$.
 \end{proposition}

\begin{proof}
The formula \eqref{doubling1} implies the square integrability since
the right hand side of \eqref{doubling1} is absolutely convergent by
Weil's convergence criterion~\eqref{Weilcrit}. We have
\begin{align*}
\|\Lambda(f)\|_2^2 =
 \int_{X} &\left( \int_{\G(N) \back \h} f(\tau_1)
\overline{\theta(\tau_1,\varphi_{q,\ell}\otimes \varphi_f)}
v_1^{\kappa} d\mu(\tau_1) \right) \\
&\quad\wedge \overline{
 \left( \int_{\G(N) \back \h} f(\tau_2) \overline{\theta(\tau_2,\ast \varphi_{q,\ell}\otimes \varphi_f)} v_2^{\kappa} d\mu(\tau_2)  \right)}.
\end{align*}
Interchanging the integration, we obtain
\[
 \int\int f(\tau_1) \overline{f(\tau_2)}
 \overline{ \left(\int_{X}  \theta(\tau_1,\varphi_{q,\ell}\otimes\varphi_f) \wedge \overline{\theta(\tau_2,\ast \varphi_{q,\ell}\otimes \varphi_f)}\right)}
(v_1v_2)^{\kappa} d\mu(\tau_1)d\mu(\tau_2).
\]
Since $\varphi_{q,\ell}$ is real valued, we easily see by the
explicit formulas of the Weil representation that
\[
\overline{\theta(\tau_2,\ast \varphi_{q,\ell}\otimes\varphi_f)} =
\theta(-\bar{\tau}_2, \ast \varphi_{q,\ell}\otimes\varphi_f)
\]
and therefore
\[
 \theta(\tau_1,\varphi_{q,\ell}\otimes\varphi_f) \wedge \overline{\theta(\tau_2,\ast
 \varphi_{q,\ell}\otimes\varphi_f)}
= \theta(\tau_1,-\bar{\tau}_2,z,\phi_{q,\ell}\otimes\phi_f) \mu
\]
by \eqref{phiKMdef}. This implies the assertion.
\end{proof}

\begin{remark}
For signature $(p,2)$, the lift $\Lambda(f)$ is actually always
square integrable, see \cite{Br2,BF}. We expect this to be true for
other signatures as well even if Weil's convergence criterion does
not hold. In that case, one would need to regularize the theta
integral $\tilde{I}$ as in \cite{KR3}.
\end{remark}

Note that the Schwartz function $\xi$ introduced by \eqref{xi:def} is
$K_{\infty}^{H}$-invariant. We can therefore consider $\xi \in
[\mathcal{S}(V^2) \otimes C^{\infty}(D)]^{H(\R)}$ by setting
\[
\xi(\x,z) = \xi(h^{-1}_{\infty}\x)
\]
with $h_{\infty} \in H(\R)$ such that $h_{\infty} z_0=z$. In
particular, $ \xi(\x,z_0)=\xi(\x)$.

\begin{proposition}\label{keyprop}
Define $\theta(\tau_1,\tau_2,z,\xi\otimes\phi_f)$ and
$\tilde{I}(\tau_1,\tau_2,\xi\otimes\phi_f)$ in the same way
as for $\phi_{q,\ell}$ in \eqref{thetadef1}, \eqref{intdef1}. Then
\begin{equation*}
\|\Lambda(f)\|_2^2 = \left(f(\tau_1) \otimes \overline{f(\tau_2)},
\tilde{I}(\tau_1,-\bar{\tau}_2,\xi\otimes\phi_f)\right).
\end{equation*}

\end{proposition}

\begin{proof}
By Proposition~\ref{decomp} and Proposition~\ref{doubling0}, we see
(omitting $\phi_f$ from the notation)
\begin{align*}
\|\Lambda(f)\|_2^2 &=   \left(f(\tau_1) \otimes \overline{f(\tau_2)}, \tilde{I}(\tau_1,-\bar{\tau}_2,\phi_{q,\ell})\right) \\
&=  \left(f(\tau_1) \otimes \overline{f(\tau_2)}, \tilde{I}(\tau_1,-\bar{\tau}_2,\xi)\right) + \left(f(\tau_1) \otimes \overline{f(\tau_2)}, \tilde{I}(\tau_1,-\bar{\tau}_2,\omega(R_1)\omega(R_2)\psi)\right)\\
&=  \left(f(\tau_1) \otimes \overline{f(\tau_2)},
\tilde{I}(\tau_1,-\bar{\tau}_2,\xi)\right) + \left(f(\tau_1) \otimes
\overline{f(\tau_2)},
R_1R_2\tilde{I}(\tau_1,-\bar{\tau}_2,\psi)\right).
\end{align*}
By the adjointness of the Maass lowering and raising operators with
respect to the Petersson scalar product, the latter summand
vanishes.
\end{proof}

\begin{corollary}
Let $p=1$ and $q+\ell>1$. Then $\Lambda$ vanishes identically.
\end{corollary}

\begin{proof}
This is obvious from Proposition~\ref{keyprop} and $\xi=0$
(Lemma~\ref{xi:vanishing}).
\end{proof}

\begin{remark}
We could have defined the lift $\Lambda$ of $f$ by using the
Schwartz form $\varphi_{q,[\ell]}$ instead of the form
$\varphi_{q,\ell}$. Using Lemma~\ref{Brauer} we see by the argument
of the proof of Proposition~\ref{keyprop} that the $L^2$-norms
$\|\Lambda(f)\|$ coincide.
\end{remark}

We want to relate the integral
$\tilde{I}(\tau_1,\tau_2,\xi\otimes\phi_f)$ to the pullback
of a genus $2$ Eisenstein series via the Siegel-Weil formula.
We first need to relate the integral over the locally
symmetric space $X$ to an integral over $H(\Q) \back H(\A)$.
We do this following \cite{KBforms}, pp.~332. First we define
the theta series associated to $\xi$ more generally for $g \in
G(\A)$ and $h=(h_{\infty}h_f) \in H(\A)$ by
\[
\theta(g,h,\xi\otimes\phi_f) = \sum_{\x \in V^2(\Q)} \omega(g)
\xi(h^{-1}_{\infty}\x,z_0)\phi_{f}(h^{-1}_{f}\x),
\]
where $z_0$ is the base point of $D$. %Note that for  $z \in D$ and $h_{\infty}(z) \in H_0(\R)$ such that $z = h_{\infty}(z) z_0$, we have $\phi_{\kappa}(h^{-1}_{\infty}(z)\x,z_0) = \phi_{\kappa}(\x,z)$.
Note that
\[
\theta(\tau_1,\tau_2,z,\xi\otimes \phi_f) = (v_1v_2)^{-\kappa/2}
\theta(\iota_0(g_{\tau_1},g_{\tau_2}),h_{\infty},\xi\otimes\phi_f)
\]
with $h_{\infty} \in H(\R)$ such that $z=h_{\infty}z_0$.

We normalize the Haar measure on $H(\A)$ such that $H(\Q) \back
H(\A)$ has volume $1$. Moreover, we normalize the Haar
measure $dh_{\infty}$ on $H(\R)$ such that
\[
\int_{D} f(z) \,\mu = \int_{H(\R)} f(h_{\infty}z_0) \,dh_{\infty}
\]
for compactly supported functions $f$ on $D$. Here $\mu$ denotes the
measure on $D$ induced by the invariant Riemann metric normalized as in
\eqref{metric-normalization}.  This gives rise to a factorization $dh=
dh_{\infty} \times dh_f$.

\begin{remark}
\label{rem:vol}
One has, see also \cite{KBforms} Remark~4.18, that
\[
\vol(K_f^H) = \frac1{\vol(X,\mu)},
\]
where $\vol(X,\mu)$ denotes the volume of $X$ with respect to the
volume form $\mu$ on $D$.
\end{remark}

\begin{proposition}
\label{prop:norm}
We have
\[
\frac{1}{\vol(X,\mu)} \tilde{I}(\tau_1,\tau_2,\xi\otimes\phi_f) =(v_1v_2)^{-\kappa/2}
\int_{H(\Q) \back
H(\A)}\theta(\iota_0(g_{\tau_1},g_{\tau_2}),h,\xi\otimes\phi_f) dh,
\]
were $dh$ is the invariant measure on $H(\A)$ such that $H(\Q) \back
H(\A)$ has volume $1$.
\end{proposition}

\begin{proof}
We use the above normalizations of the Haar measures on $H(\R)$ and $H(\A_f)$.
By means of  \eqref{strong-approx} we obtain
\begin{align*}
\int_{H(\Q) \back
H(\A)}&\theta(\iota_0(g_{\tau_1},g_{\tau_2}),h,\xi\otimes\phi_f) dh
\\&= \sum_j \int_{H(\Q) \back H(\Q) H(\R)_0 h_j K_f^H
h_j^{-1}}\theta(\iota_0(g_{\tau_1},g_{\tau_2}),hh_j,\xi\otimes\phi_f)
dh \\ &= \vol(K_f^H) \sum_j \int_{\G_j \back
H(\R)_0}\theta(\iota_0(g_{\tau_1},g_{\tau_2}),h_{\infty}
h_j,\xi\otimes\phi_f) dh_{\infty} \\
&=(v_1v_2)^{\kappa/2}\vol(K_f^H)\tilde{I}(\tau_1,\tau_2,\xi\otimes\phi_f).
\end{align*}
The assertion now follows using Remark \ref{rem:vol}.
\end{proof}

%We now assume that $\xi \otimes \phi_f$ is "even", i.e.,
%\[ \theta(\iota_0(g_{\tau_1},g_{\tau_2}),\eps
%h_{\infty},\xi\otimes\varphi_f) =
%\theta(\iota_0(g_{\tau_1},g_{\tau_2}),h_{\infty},\xi\otimes\varphi_f)
%\]
%for any $\eps \in H(\R) /H(\R)_0$. This is for example the case if
%$q+\ell$ is even and $\varphi_f$ corresponds to the characteristic
%function of $L$.
%From now on, we assume that
%\[
%\text{$L$ is unimodular,}
%\]
%and therefore in particular
%\[
%\text{$q$ is even.}
%\]

\begin{proposition}
\label{prop:l21}
Let $\Xi(s) \otimes \Phi_f(s)$ be the section associated to
$\xi\otimes\phi_f$ via \eqref{Section} and let $s_0 = (m-3)/2$. Then
\[
\frac{1}{\vol(X,\mu)} \|\Lambda(f)\|_2^2 = (v_1v_2)^{-\kappa/2}
\left(f(\tau_1) \otimes \overline{f(\tau_2)},
E(\iota_0(g_{\tau_1},g_{-\bar{\tau}_2}),s_0, \Xi \otimes \Phi_f)\right).
\]
\end{proposition}

\begin{proof}
Using Proposition \ref{prop:norm} and the Siegel-Weil formula,  Theorem~\ref{Siegel-Weil},  we find
\[
\frac{1}{\vol(X,\mu)}
\tilde{I}(\tau_1,\tau_2,\xi\otimes\phi_f) =(v_1v_2)^{-\kappa/2}
E(\iota_0(g_{\tau_1},g_{\tau_2}),s_0, \Xi \otimes
\Phi_f).
\]
Now the assertion follows from Proposition~\ref{keyprop}.
\end{proof}

\begin{corollary}\label{cor:q+ell-even}
Assume that $q+\ell$ is even and $p>1$. Let $\Phi_\infty^\kappa(s)$
be the standard section defined by \eqref{eq:phikappa}, and let
$\Phi_f(s)$ be the section associated to $\phi_f$ via
\eqref{Section}. Then
\[
\frac{1}{\vol(X,\mu)} \|\Lambda(f)\|_2^2 = C(s_0)
(v_1v_2)^{-\kappa/2} \left(f(\tau_1) \otimes \overline{f(\tau_2)},
E(\iota_0(g_{\tau_1},g_{-\bar{\tau}_2}),s_0,
\Phi_\infty^\kappa\otimes \Phi_f)\right),
\]
where $C(s_0)$ is the nonzero constant in
Proposition~\ref{Xi:Indent}.
%In particular,
%\[
%C(s_0) \ne 0 \Longleftrightarrow p>1.
%\]
\end{corollary}

\begin{proof}
We have $\Xi(g,s) = C(s) \Phi_{\infty}^{\kappa}(g,s)$ by
Proposition~\ref{Xi:Indent}. Hence the Corollary immediately follows
from Proposition \ref{prop:l21}.
\end{proof}

Suppose that $f$ is an eigenform of level $N$ and let $S$ denote the
set of primes dividing $N$ together with $\infty$.  Then the
doubling method \cite{PS-R, Ra, Boe, Ga} expresses a convolution
integral as on the right hand side above as a product of the
standard $L$-function $L^S(s_0+\frac{1}{2},f)$ with the Euler
factors corresponding to $p\in S$ omitted times a product of ``bad''
local factors corresponding to the primes in $S$.  If $m>4$ then
$s_0+\frac{1}{2}$ lies in the region of convergence of the Euler
product of $L^S(s,f)$. Hence the $L$-value does not vanish.
Therefore the lift $\Lambda(f)$ vanishes precisely if at least one
of the ``bad'' local factors vanishes.  By the analysis of the
present paper we determine the local factor at infinity.

We now specialize to the case when the lattice $L$ is even and {\em
  unimodular}.  Then $\varphi_f$ corresponds to the characteristic
function of $L$ and $\Phi_f(s) =1$.  The level of $L$ is $N=1$, so
that $\infty$ is the only ``bad'' place. By the above analysis we
obtain a very explicit formula for $\|\Lambda(f)\|_2^2$ as we shall now
explain.

In this case $\theta(\tau,z,\varphi_{q,\ell})$ is a modular form of weight
 $\kappa = m/2+\ell$ for $\SL_2(\Z)$
and vanishes unless $q+\ell$ is even,
which we
assume from now on as well.
Then $\kappa$ is even, because $8\mid p-q$.
By Corollary \ref{cor:q+ell-even} and \eqref{classicaleis} we have
\begin{align}
\label{eq:doublingclassical} \frac{1}{\vol(X,\mu)}
\|\Lambda(f)\|_2^2 = C(s_0) \left(f(\tau_1) \otimes
\overline{f(\tau_2)},
E^{(2)}_{\kappa}(\tau_1,-\bar{\tau}_2,-\ell/2)\right),
\end{align}
where $E^{(2)}_{\kappa}(\tau_1,\tau_2,s)$ is the pullback of the classical
genus 2 Siegel Eisenstein series $E_\kappa^{(2)}(\tau,s)$ (see
\eqref{classicaleis2}) to the diagonal.

% \texttt{ vielleicht jetzt ein allgemeines statement, dass da nun
% bis auf endlich viele Faktoren ein Wert der Standard L-funktion
% rauskommt}

%\texttt{hier oder in die Einleitung die Erlaueterung, dass man bei
%Anwendungen nicht an den endlichen Stellen rumspielen darf; unsere
%Rechnungen zeigen auf jeden Fall, dass bei unendlich nichts
%passiert}

We recall the definition of the standard $L$-function of a Hecke
eigenform  $f\in S_\kappa(\Gamma(1))$.
We use the normalization of \cite{An},
\cite{Boe},  \cite{La}.
We denote the Fourier
coefficients of $f$ by $c(n)$ and assume that $f$ is normalized, i.e.,
$c(1)=1$.
Let $p$ be a prime. The Satake parameters $\alpha_{0,p}, \alpha_{1,p}$
of $f$ at $p$ are defined by
the factorization of the Hecke polynomial
\begin{align}
%Q_p(X)=
(1-c(p)X+p^{\kappa-1} X^2)=(1-\alpha_{0,p} X)(1-\alpha_{0,p}\alpha_{1,p} X).
\end{align}
Hence
\begin{align*}
\alpha_{0,p}^2\alpha_{1,p}&=p^{\kappa-1},  &\alpha_{0,p}(1+\alpha_{1,p})&=c(p).
\end{align*}
According to Deligne's theorem, formerly the Ramanujan-Petersson conjecture, we
have $|\alpha_{1,p}|=1$.
The standard $L$-function of $f$ is defined
by the Euler product
\begin{align}
\label{eq:standard}
D_f(s)=\prod_p \big[ (1-p^{-s})(1-\alpha_{1,p}^{-1} p^{-s})(1-\alpha_{1,p} p^{-s})\big]^{-1}.
\end{align}
It converges for $\Re(s)>1$.
The corresponding completed $L$-function
\begin{align}
\Psi_f(s)=\pi^{-\frac{3s}{2}}\Gamma\left(\frac{s+1}{2}\right)\Gamma\left(\frac{s+\kappa-1}{2}\right)\Gamma\left(\frac{s+\kappa}{2}\right) D_f(s)
\end{align}
has a  meromorphic continuation to $\C$ and satisfies the functional equation
\begin{align}
\Psi_f(s)=\Psi_f(1-s)
\end{align}
(see e.g.~\cite{Boe}, \cite{Sh}).
%
\begin{comment}
For completeness we also recall the
different normalization used in \cite{Sh} and the relation to the
Fourier coefficients of $f$.
If we put
\begin{align*}
\alpha_p&=\alpha_{0,p},  &\beta_p&= \alpha_{0,p}\alpha_{1,p},
\end{align*}
%Then we have
%\begin{align*}
%\alpha_p\beta_p&=p^{\kappa-1},  &\alpha_p+\beta_p&= c(p).
%\end{align*}
we may rewrite $D_f(s)$ as follows:
\[
D_f(s)=\prod_p \big[ (1-\alpha_p\beta_p p^{1-\kappa-s})(1-\alpha_{p}^{2} p^{1-\kappa-s})(1-\beta_p^2 p^{1-\kappa-s} )\big]^{-1}.
\]
Comparing this with the definition of the $L$-function $D(s,f,\chi_0)$ in
\cite{Sh}, we see that $D_f(s)=D(s+\kappa-1,f,\chi_0)$.
\end{comment}
It is well known  (see \cite{Sh}, Introduction, \cite{Za}) that $D_f(s)$ can be
interpreted as the Rankin $L$-series
\begin{align*}
D_f(s) &=\zeta(2s)\sum_{n=1}^\infty c(n^2) n^{-s-\kappa+1}
%\\&
=\frac{\zeta(2s)}{\zeta(s)}\sum_{n=1}^\infty c(n)^2 n^{-s-\kappa+1}.
\end{align*}

\begin{theorem}
\label{thm:main} Assume that $m>3+r$ so that Weil's convergence
criterion~\eqref{Weilcrit} in genus $2$ holds. Furthermore, assume
that $q+\ell$ is even and that $L$ is even unimodular. Let $f\in
S_\kappa(\Gamma(1))$ be a Hecke eigenform, and write
$\|f\|^2_2=(f,f)$ for its Petersson norm normalized as in
\eqref{petnorm}.
%and write $D_f(s)$
%for the standard $L$-function of $f$ normalized as in \cite{Boe}.
We have
\begin{align*}
\frac{1}{\vol(X,\mu)} \cdot \frac{\|\Lambda(f)\|_2^2}{\|f\|^2_2} =
C(s_0) \mu(1,\kappa,-\ell/2)\frac{D_f(m/2-1)}{\zeta(m/2)\zeta(m-2)},
\end{align*}
where
\begin{align*}
\mu(1,\kappa,-\ell/2)=2^{3-m/2} (-1)^{\kappa/2} \pi \frac{\Gamma(m/2+\ell/2-1)}{\Gamma(m/2+\ell/2)}.
\end{align*}
\end{theorem}

\begin{proof}
The statement follows from \eqref{eq:doublingclassical} by means of
\cite{Boe}, identities (14) and (22).
\end{proof}

\begin{remark}
  By the same argument it is easily seen that
  $(\Lambda(f),\Lambda(g))=0$ for two different normalized Hecke
  eigenforms $f$ and $g$.
\end{remark}

\begin{corollary}
  Assume that $m>\max(4,3+r)$, $p>1$, $q+\ell$ even, and that $L$ is even unimodular.  Then
  the theta lift $\Lambda: S_\kappa(\Gamma(1))\to
  \mathcal{Z}^q(X,\widetilde{\Sym}^{\ell}(V))$ is injective.
\end{corollary}

\begin{proof}
  This follows from Theorem \ref{thm:main}, Proposition \ref{Xi:Indent}, and the convergence of the
  Euler-product for $D_f(m/2-1)$ in this case.
\end{proof}

\begin{comment}

 We now specialize to the case when $L$ is unimodular. Then $\varphi_f$ corresponds to the
 characteristic function of $L$ and $\Phi_f(s) =1$. Moreover,
 $\theta(\tau,z,\varphi_{q,\ell})$ is a modular form of weight
 $\kappa = m/2+\ell$ for $\SL_2(\Z)$ and vanishes unless $q+\ell$ is even, which we
 assume from now on as well. By \eqref{xiweight} we have $\Xi(g,s) = \Phi_{\infty}^{\kappa}(g,s)$.

 The classical Siegel Eisenstein series of weight $\kappa =
 m/2+\ell$ for $\Sp(2,\Z)$ is given by
 \[
 E^{(2)}_{\kappa}(\tau,s) = \sum_{\g \in \G_{\infty} \back \Sp(2,\Z)}
\left(\frac{\det(v)}{|\det(c\tau+d)|^2}\right)^s
\det(c\tau+d)^{-\kappa},
\]
where $\g = \kzxz{a}{b}{c}{d}$. Here $\tau = u+iv \in \h_2$. Thus
\[
\det(v)^{-\kappa/2} E(g_{\tau},s,\Phi^{\kappa}_{\infty}(s)\otimes
\Phi_f) =E^{(2)}_{\kappa}\left(\tau,\tfrac{s+3/2-\kappa}{2}\right)
\]
and
\[
\|\Lambda(f)\|_2^2 =  \frac1{vol(K_f^H)} \left(f(\tau_1) \otimes
\overline{f(\tau_2)},E^{(2)}_{\kappa}(\tau_1,\tau_2,-\ell/2)\right).
\]
where $E^{(2)}_{\kappa}(\tau_1,\tau_2,s)$ is pullback of $E^{(2)}$
to the diagonal.

\texttt{Wie lautet denn hier nun die Funktionalgleichung?}
\end{comment}


\begin{thebibliography}{99}



%\bibitem{AbSt} M. Abramowitz and I. Stegun, {\em Pocketbook of Mathematical Functions}, Verlag Harri Deutsch (1984).


\bibitem{An} A. N. Andrianov, \emph{The multiplicative arithmetic of Siegel
modular forms}, Russian Math. Surveys {\bf 34} (1979), 75--148.

\bibitem{Bo1} R. E. Borcherds, \emph{Automorphic forms on $\Orth_{s+2,2}(\R)$
and infinite products}, Invent. Math. {\bf 120} (1995), 161--213.

\bibitem{Bo2}
R. E. Borcherds, \emph{Automorphic forms with singularities on
Grassmannians}, Inv. Math. \textbf{132} (1998), 491--562.

\bibitem{Bo3}
R. Borcherds, \emph{The Gross-Kohnen-Zagier theorem in higher
dimensions}, Duke Math. J. \textbf{97} (1999), 219--233. Correction
in: Duke Math J. \textbf{105} No. 1 p.183--184.

%\bibitem{Bost}
%J.-B. Bost,
%\emph{Potential theory and Lefschetz theorems for arithmetic surfaces}, Ann. Sci \'Ecole Norm. Sup. \textbf{32} (1999), 241-312.




%\bibitem{Br1}
%J. Bruinier,
%\emph{Borcherds products and Chern classes of Hirzebruch-Zagier divisors},
%Inv. Math. \textbf{138} (1999), 51-83.

\bibitem{Boe} S. B\"ocherer, \emph{\"Uber die Funktionalgleichung
  automorpher $L$-Funktionen zur Siegelschen Modulgruppe}, J. Reine
  Angew. Math.  362 (1985), 146--168.


\bibitem{Br2} J. Bruinier, \emph{Borcherds products on $\Orth(2,l)$
    and Chern classes of Heegner divisors}, Springer Lecture Notes in
  Mathematics {\bf 1780}, Springer-Verlag (2002).

%\bibitem{BBK}
%J. Bruinier, J. Burgos, and U. K\"uhn,
%\emph{Borcherds products and arithmetic intersection theory on Hilbert modular surfaces}, preprint.


\bibitem{BrFr} {\em J. H. Bruinier and E. Freitag}, Local Borcherds products,
Annales de l'Institut Fourier {\bf 51.1} (2001), 1--26.


\bibitem{BF}
J. Bruinier and J. Funke, \emph{On two geometric theta lifts}, Duke
Math J. \textbf{125} (2004), 45-90.



%\bibitem{BKK}
%J. Burgos, J. Kramer, and U. K\"uhn,
%\emph{Cohomological Arithmetic Chow groups}, preprint (2003).
%arXiv.org/math.AG/0404122


%\bibitem{E} A. Erd\'elyi, W. Magnus, F. Oberhettinger and F. G. Tricomi, {\em Tables of Integral Transforms, vol.~I}, McGraw-Hill (1954).

%\bibitem{EZ} Eichler and D. Zagier, {\em M. The Theory of Jacobi Forms}, Progress in Math. {\bf 55}, Birkh\"auser (1985).


\bibitem{Fr}
E. Freitag, \emph{Stabile Modulformen}, Math. Ann. {\bf 230} (1977), 197--211.


\bibitem{FH}
W. Fulton and J. Harris, \emph{Representation Theory, A First
Course}, Graduate Texts in Mathematics {\bf 129}, Springer, 1991.



%\bibitem{FCompo}
%J. Funke,
%\emph{Heegner divisors and nonholomorphic modular forms}, Compositio Math. \textbf{133} (2002), 289-321.

\bibitem{FM}
J. Funke and J. Millson, \emph{Cycles with local coefficients for
orthogonal groups and vector-valued Siegel modular forms}, to appear
in American J. Math. (2006).

\bibitem{Ga} P. Garrett, \emph{Pullbacks of Eisenstein series;
  Applications}. In: Automorphic forms of several variables, Taniguchi
  Symposium, Katata, 1983, Birh\"auser (1984).

%\bibitem{GH}
%P. Griffiths and J. Harris, \emph{Principles of Algebraic Geometry}, Wiley (1978).

\bibitem{HM} {\em J. Harvey and G. Moore}, Algebras, BPS states, and strings, Nuclear Phys. B {\bf 463} (1996), no. 2-3, 315--368.

%\bibitem{He} B. E. Heim, \emph{ Pullbacks of Eisenstein series, Hecke-Jacobi theory and automorphic $L$-functions}, Proc. Sympos. Pure Math. {\bf 66}, Part 2, Amer. Math. Soc., Providence, RI, 201--238 (1999).



%\bibitem{HZ}
%F. Hirzebruch and D. Zagier,
%\emph{Intersection numbers of curves on Hilbert modular surfaces and modular forms of Nebentypus},
%Inv. Math. \textbf{36} (1976), 57-113.

%\bibitem{KS}
%S. Katok and P. Sarnak, \emph{Heegner points, cycles and Maass forms}, Israel J. Math. \textbf{84} (1993), 193-227.

%\bibitem{Kim1}
%C. H. Kim, \emph{Borcherds products associated with certain Thompson series},
%Compositio Math. \textbf{140} (2004), 541-551.

%\bibitem{Kim2}
%C. H. Kim, \emph{Traces of singular moduli and Borcherds products},
%preprint (2003).

%\bibitem{Kuehn}
%U. K\"uhn,
%\emph{Generalized arithmetic intersection numbers},
%J. reine angew Math. 534 (2001), 209-236.



%\bibitem{Franke}
%H.-G. Franke,
%\emph{Kurven in Hilbertschen Modulfl\"achen und Humbertsche Fl\"achen},
%Bonner math. Schriften, vol. 114, 1978.


%\bibitem{Hausmann}
%W. Hausmann,
%\emph{Kurven auf Hilbertschen Modulfl\"achen},
%Bonner math. Schriften, vol. 123, 1980.

\bibitem{Li} J.-S. Li, %Jian-Shu Li
\emph{Nonvanishing theorems for the cohomology of certain arithmetic quotients},  J. Reine Angew. Math.  {\bf 428}  (1992), 177--217.

%\bibitem{Kubota}
%T. Kubota, \emph{Elementary Theory of Eisenstein Series}, Halsted
%Press, (1973).



%\bibitem{KShintani}
%S. Kudla, \emph{ On the integrals of certain singular
%theta-functions}, J. Fac. Sci. Univ. Tokyo \textbf{28} (1982),
%439-465.


\bibitem{Ku}
S. Kudla, \emph{Some extensions of the Siegel-Weil formula},
unpublished manuscript (1992). Available at {\tt
www.math.umd.edu/\~{}ssk} .



%\bibitem{KDuke}
%S. Kudla,
%\emph{Algebraic cycles on Shimura varieties of orthogonal type},
%Duke Math. J. \textbf{86} (1997), 39-78.

%\bibitem{KAnn}
%S. Kudla,
%\emph{Central derivatives of Eisenstein series and height pairings},
%Ann. of Math. \textbf{146} (1997), 545-646.

\bibitem{KBforms}
S. Kudla, \emph{Integrals of Borcherds forms}, Compositio Math.
\textbf{137} (2003), 293-349.


%\bibitem{Kmsri}
%S. Kudla,
%\emph{Special cycles and derivatives of Eisenstein series},
%MSRI proceeding on Heegner points (to appear).



%\bibitem{KBourbaki}
%S. Kudla, \emph{Derivatives of Eisenstein series and generating
%functions for arithmetic cycles}, S\'eminaire Bourbaki,
%\textbf{876} (2000).

\bibitem{KMI}
S. Kudla and J. Millson, \emph{The Theta Correspondence and Harmonic
Forms I}, Math. Ann. \textbf{274} (1986), 353-378.


%\bibitem{KMII}
%S. Kudla and J. Millson,
%\emph{The Theta Correspondence and Harmonic Forms II},
%Math. Ann. \textbf{277} (1987), 267-314.

%\bibitem{KMCan}
%S. Kudla and J. Millson, \emph{Tubes, cohomology with growth
%conditions and an application to the theta correspondence}, Canad.
%J. Math.  \textbf{40}  (1988), 1-37.

\bibitem{KM90}
S. Kudla and J. Millson,
\emph{Intersection numbers of cycles on locally symmetric spaces and Fourier coefficients of holomorphic modular forms in several complex variables},
IHES Pub. \textbf{71} (1990), 121-172.

\bibitem{KR1}
S. Kudla and S. Rallis, \emph{On the Weil-Siegel formula},  J. Reine
Angew. Math. \textbf{387} (1988), 1-68.

\bibitem{KR2}
S. Kudla and S. Rallis, \emph{On the Weil-Siegel formula II},  J.
Reine Angew. Math. \textbf{391} (1988), 65-84.

\bibitem{KR}
S. Kudla and S. Rallis, \emph{Degenerate principal series and
invariant distributions}, Israel J. Math. \textbf{69} (1990), 25-45.

\bibitem{KR3}
S. Kudla and S. Rallis, \emph{A regularized Weil-Siegel formula: the
first term identity},  Annals of Math. \textbf{140} (1994), 1-80.



%\bibitem{KRY}
%S. Kudla, M. Rapoport and T. Yang,
%\emph{On the derivative of an Eisenstein series of weight one},
%Int. Math. Res. Not. \textbf{7} (1999), 347-385.

%\bibitem{KRY}
%S. Kudla, M. Rapoport and T. Yang,
%\emph{Derivatives of Eisenstein series and Faltings heights}, Compositio Math. (to appear).

\bibitem{La} R. P. Langlands, \emph{Problems in the theory of
    automorphic forms}, Lecture Notes in Math. {\bf 170} (1970),
  18--86, Springer-Verlag.




%\bibitem{Le}
%N.N. Lebedev,
%{Special functions and their applications},
%Dover, 1972.


%\bibitem{LV}
%G. Lion and M. Vergne,
%\emph{The Weil representation, Maslov index and theta series},
%Progress in Math., vol. 6, Birkh\"auser 1980

%\bibitem{Maass}
%H. Maass, \emph{\"Uber die r\"aumliche Verteilung der Punkte in Gittern mit indefiniter Metrik}, Math. Ann. \textbf{138} (1959), 287-315.


%\bibitem{Miyake}
%T. Miyake,
%\emph{Modular forms},
%Springer, 1989.




%\bibitem{Oda}
%T. Oda,
%\emph{On modular forms associated with indefinite quadratic forms of signature $(2,n-2)$},
%Math. Annalen \textbf{231} (1977), 97-144.

\bibitem{PS-R} I. Piatetski-Shapiro, S. Rallis, $L$-functions for classical groups.
Lecture Notes in Mathematics {\bf 1254},
Springer-Verlag, Berlin (1987).

\bibitem{Ra} S. Rallis, \emph{Injectivity properties of liftings associated to Weil representations},  Compositio Math.  {\bf 52}  (1984), 139--169.

%\bibitem{RS}
%S. Rallis and G. Schiffmann,
%\emph{On a relation between $\widetilde{SL_2}$ cusp forms and cusp forms on tube domains associated to orthogonal groups},
%Trans. AMS \textbf{263} (1981), no.1, 1-58.



\bibitem{Sh}
G. Shimura, \emph{On the holomorphicity of certain Dirichlet series},
Proc. London Math. Soc. \textbf{31} (1975), 79--98.





%\bibitem{Soule}
%C. Soul\'e at al., \emph{Lectures on Arakelov Geometry}, Cambridge Studies in Advanced Mathematics \textbf{33}, Cambridge University Press (1992).


%\bibitem{vG}
%G. van der Geer,
%\emph{Hilbert modular surfaces},
%Ergebnisse der Math. und ihrer Grenzgebiete (3), vol. 16, Springer, 1988.

%\bibitem{Wang85}
%S.P. Wang,
%\emph{Correspondence of modular forms to cycles associated to $O(p,q)$},
%J. Diff. Geom. \textbf{18} (1985), 151-223.

%\bibitem{Yang}
%T. Yang,
%\emph{Faltings heights and the derivatives of Zagier's Eisenstein series}, MSRI proceeding on Heegner points (to appear).

%\bibitem{Zagier}
%D. Zagier,
%\emph{Nombres de classes et formes modulaires de poids $3/2$},
%C. R. Acad. Sci. Paris S\'er. A-B \textbf{281} (1975), 883-886.

%\bibitem{ZagierTr}
%D. Zagier, \emph{Traces of singular moduli}, in: Motives, Polylogarithms and Hodge Theory (Part I), Eds.: F. Bogomolov and L. Katzarkov, International Press, Somerville (2002).

\bibitem{Wald}
J.-L. Waldspurger, \emph{Engendrement par des s\'eries t\^eta de certains espaces de formes modulaires}, Invent. Math. \textbf{50} (1979), 135-168.


\bibitem{Weil}
A. Weil, \emph{Sur la formule de Siegel dans la th\'eorie des
groupes classiques}, Acta Math. \textbf{113} (1965) 1-87.

\bibitem{Za}
D. Zagier, \emph{Modular forms whose Fourier coefficients involve zeta-functions of quadratic fields}. In: Modular Functions of One Variable VI, Lecture Notes in Math. {\bf 627}, Springer-Verlag (1977), 105-169.



\end{thebibliography}
\end{document}